%% file: ms.tex
\documentclass[8pt,twocolumn]{extarticle}

\usepackage{algorithm}
\usepackage{algorithmic}
\usepackage{amssymb}
\usepackage[toc, page, title, titletoc]{appendix}
\usepackage[left=1.5cm, right=1.5cm]{geometry}
\usepackage{graphicx}
\usepackage[hyphens]{url}
\usepackage[hidelinks]{hyperref}
\usepackage[hyphenbreaks]{breakurl}
\usepackage[cp1252]{inputenc}
\usepackage{listings}
\usepackage{mathtools}
\usepackage[square]{natbib}
\usepackage{xcolor}

\DeclarePairedDelimiter\floor{\lfloor}{\rfloor}

\setcitestyle{square,aysep={}}
 
\begin{document}

\title{Generalising the Fast Reciprocal Square Root Algorithm}
\author{Mike Day \\ \href{mailto:mike.day@cantab.net}{mike.day@cantab.net}}
\date{}
\maketitle

\input{abstract}

\input{introduction}

\input{related_work}
\input{standard_algorithm}
\section{Analysis and Generalisation} \label{sec:analysis-and-generalisation}
\input{coarse_approximation}
\input{refined_approximation}
\input{finding_extrema}
\input{optimal_c}
\input{minimising_rho}
\input{main_result}
\input{use_of_monics}
\input{linear_minimax}
\input{multiple_iterations}
\input{implementation}
\input{results}
\input{acknowledgements}

\input{ms.bbl}
\appendix
\input{appendices}

\end{document}

%% file: abstract.tex
\begin{abstract}

The Fast Reciprocal Square Root Algorithm is a well-established approximation technique consisting of two stages: first, a coarse approximation is obtained by manipulating the bit pattern of the floating point argument using integer instructions, and second, the coarse result is refined through one or more steps, traditionally using Newtonian iteration but alternatively using improved expressions with carefully chosen numerical constants found by other authors. The algorithm was widely used before microprocessors carried built-in hardware support for computing reciprocal square roots. At the time of writing, however, there is in general no hardware acceleration for computing other fixed fractional powers. This paper generalises the algorithm to cater to all rational powers, and to support any polynomial degree(s) in the refinement step(s), and under the assumption of unlimited floating point precision provides a procedure which automatically constructs provably optimal constants in all of these cases. It is also shown that, under certain assumptions, the use of monic refinement polynomials yields results which are much better placed with respect to the cost/accuracy tradeoff than those obtained using general polynomials. Further extensions are also analysed, and several new best approximations are given.

\end{abstract}

%% file: introduction.tex
\section{Introduction}

There is a well-known algorithm for calculating a fast approximation to the reciprocal square root function using single-precision floating point arithmetic, the key to which is to treat the bits of a floating point number as if they represent an integer, and to manipulate them using integer instructions. (For an overview of this class of techniques, see \citet{blinn1997}). When the result is reinterpreted as a floating point number, it holds a coarse approximation to the reciprocal square root which can then be refined using Newtonian iteration. The algorithm was valuable in the 1990s, because it enabled rapid approximate normalisation of vectors for use in lighting calculations in graphics engines for video games.

Subsequent work (notably by \citet{lomont2003}, \citet{pizer2008}, \citet{kadlec2010}, \citet{moroz2016}, and \citet{walczyk2021}) analysed the code and found that by changing its numerical constants a roughly threefold reduction in error could be achieved, without changing the execution cost of the algorithm. However, no overarching mathematical framework was developed which would allow the same technique to be applied to the approximation of arbitrary powers of $x$, or which would automatically generate optimal coefficients for polynomials of arbitrary degree in the refinement steps.

Present-day microprocessors lessen the usefulness of the original technique since many hardware platforms provide one or more machine-level instructions to accelerate the calculation of reciprocal square roots. However, since the library function \texttt{pow()} is usually a much costlier operation, it may yet be beneficial to develop approximations which extend the technique to work for other fixed powers of $x$ (for example, those which arise in gamma correction).

%% file: related_work.tex
\section{Related Work}

An investigation by \citet{sommefeldt2006} into the history of the algorithm concluded that it was developed by Greg Walsh and Cleve Moler, who had learned about the bit-manipulation technique from an unpublished 1986 paper by William Kahan and K. C. Ng, a copy of which can be seen in the source for \texttt{fdlibm} from \citet{sun1993}. The algorithm gained widespread attention in 2005 when id Software released the source code for their game \textit{Quake III: Arena} \citep{id1998}. Their code contained a version with the ``magic constant" \texttt{0x5F3759DF}, achieving a peak relative error of $1.752339 \times 10^{-3}$. \citet{lomont2003} used analysis coupled with a numerical search to show that the optimal choice for the constant is \texttt{0x5F375A86}, lowering the error bound slightly to $1.751302 \times 10^{-3}$. \citet{moroz2016} used a purely analytical method to find the same constant Lomont had found, without requiring a numerical search. Far more significant improvements became possible using the observation that the Newtonian iteration step itself contained two further constants, which could be tuned together with the magic constant. \citet{pizer2008} thus found a trio of constants which lowered the error bound to $6.531342 \times 10^{-4}$, a $2.7$-fold improvement. Using a similar approach, \citet{kadlec2010} lowered the error bound a little further to $6.501967 \times 10^{-4}$, but still required a lengthy numerical search for the values of the constants, and no demonstration of their optimality. \citet{walczyk2021} provided an analytical method of finding a theoretically optimal set of constants, but did not perform any fine tuning to reduce evaluation error, so that Kadlec's remains the most accurate hitherto published version\footnote{When comparing results, care must be taken to ensure that the testing conditions match, since the hardware platform, the function implementation, and the instructions emitted by the compiler can all affect the outcome, even when IEEE 754 compliance is guaranteed. Our basis for comparison is described in the results section.}. The work by Walzcyk et al. also showed how to extend the optimisation technique to additional iterations, providing a dramatic improvement in accuracy over the multi-iteration version of the Quake code. A reciprocal cube root counterpart to the original algorithm was given by \citet{levin2012}. \citet{moroz2021} provided the hitherto best published constants for this algorithm, and showed how to replace the Newtonian iteration step with a quadratic Householder iteration for much greater accuracy. They proceeded to find superior coefficients for this method, but we show that their approach did not produce an optimal value for the magic constant. \citet{blinn1997} showed a corresponding coarse approximation for a general power of $x$, but did not consider optimisation of the relevant constant, or refinement of the coarse value.

%% file: standard_algorithm.tex
\section{Standard algorithm} \label{sec:standard_alg}

Much of the literature on the original technique refers to it as the ``Fast Inverse Square Root" algorithm, but since this use of the word ``inverse" is somewhat ambiguous, we will refer to it as the Fast Reciprocal Square Root (FRSR) algorithm. The Quake code for FRSR is paraphrased in Listing \ref{lst:Quake}.

\begin{minipage}{\linewidth}
\begin{lstlisting}[caption=Quake FRSR algorithm, label={lst:Quake}, language=C]
float FRSR_Quake(float x)
{
    int X = *(int *)&x;
    int Y = 0x5F3759DF - (X>>1);
    float y = *(float *)&Y;
    return y * (1.5f - 0.5f*x*y*y);
}
\end{lstlisting}
\end{minipage}

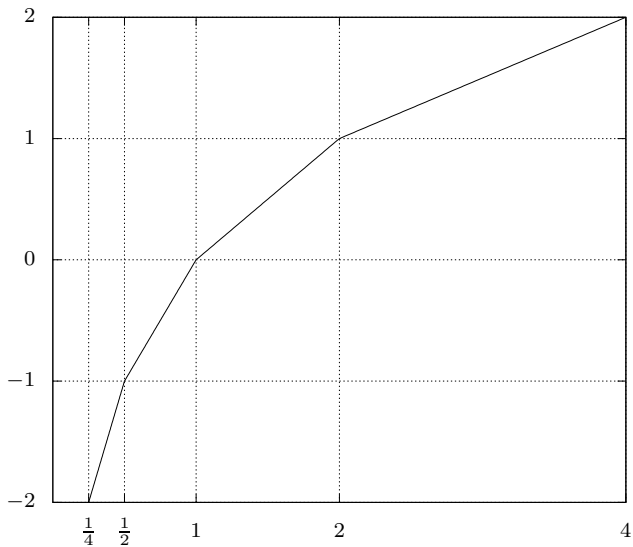
\begin{figure}
  \centering
  \input{pseudolog}
  \caption{Pseudolog function}
  \label{fig:pseudolog}
\end{figure}

It makes use of an important property of the distribution of numbers encoded in the IEEE 754 floating point format, which has been expressed by \citet{blinn1997} as follows: ``If you only deal with positive numbers, the bit pattern of a floating-point number, interpreted as an integer, gives a piecewise linear approximation to the logarithm function''.\footnote{This property breaks down for subnormal floating point arguments, so that the FRSR algorithm only applies to positive normal floats. If it is required to support subnormal arguments, a conditional branching approach such as the one described in \citet{walczyk2021} can be used.} Strictly, the mapping also introduces a scale and bias, but this fact is not relevant to the analysis since the scale is cancelled out when the inverse bit interpretation occurs, and the effect of the bias can be counteracted by suitably modifying the hexadecimal constant. Hence, we will use $L(x), \hspace{2pt} x>0$ to mean the piecewise linear function which has the value $\log_2{x}$ at points where $x$ is a a power of two, and is linear between each successive pair of such points. We will call $L(x)$ the \textit{pseudolog} of $x$ (see Figure \ref{fig:pseudolog}). A major advantage of using this approach is that it decouples the analysis from the floating point representation, and we can now treat all values as real numbers. The standard FRSR algorithm can then be expressed as shown in Algorithm \ref{alg:FRSR}.

\begin{algorithm}
\caption{Standard FRSR algorithm}
\label{alg:FRSR}
\begin{algorithmic}[1]
\STATE \textbf{function} FRSR($x$)
\STATE $X = L(x)$                                      \label{alg:FRSR:line:X}
\STATE  $Y = \frac{c}{2}-\frac{1}{2}X$   \label{alg:FRSR:line:Y}
\STATE  $y = L^{-1}(Y)$                           \label{alg:FRSR:line:y}
\STATE  $z = xy^2$                                  \label{alg:FRSR:line:z}
\STATE  \textbf{return} $ y(c_0+c_1z) $ \label{alg:FRSR:line:return}
\end{algorithmic}
\end{algorithm}

The constant $\frac{c}{2}$ of line \ref{alg:FRSR:line:Y} corresponds, via a scale and bias, to FRSR's ``magic constant'' (the reason for dividing $c$ by $2$ will become apparent), and the coefficients $c_0$ and $c_1$ represent the generalised values which researchers have used to improve upon the $\frac{3}{2}$ and $-\frac{1}{2}$ of the Newtonian iteration of Listing \ref{lst:Quake}.

%% file: pseudolog.tex
\begingroup
  \inputencoding{cp1252}%
  \fontfamily{Times-New-Roman}%
  \selectfont
  \makeatletter
  \providecommand\color[2][]{%
    \GenericError{(gnuplot) \space\space\space\@spaces}{%
      Package color not loaded in conjunction with
      terminal option `colourtext'%
    }{See the gnuplot documentation for explanation.%
    }{Either use 'blacktext' in gnuplot or load the package
      color.sty in LaTeX.}%
    \renewcommand\color[2][]{}%
  }%
  \providecommand\includegraphics[2][]{%
    \GenericError{(gnuplot) \space\space\space\@spaces}{%
      Package graphicx or graphics not loaded%
    }{See the gnuplot documentation for explanation.%
    }{The gnuplot epslatex terminal needs graphicx.sty or graphics.sty.}%
    \renewcommand\includegraphics[2][]{}%
  }%
  \providecommand\rotatebox[2]{#2}%
  \@ifundefined{ifGPcolor}{%
    \newif\ifGPcolor
    \GPcolorfalse
  }{}%
  \@ifundefined{ifGPblacktext}{%
    \newif\ifGPblacktext
    \GPblacktexttrue
  }{}%
  \let\gplgaddtomacro\g@addto@macro
  \gdef\gplbacktext{}%
  \gdef\gplfronttext{}%
  \makeatother
  \ifGPblacktext
    \def\colorrgb#1{}%
    \def\colorgray#1{}%
  \else
    \ifGPcolor
      \def\colorrgb#1{\color[rgb]{#1}}%
      \def\colorgray#1{\color[gray]{#1}}%
      \expandafter\def\csname LTw\endcsname{\color{white}}%
      \expandafter\def\csname LTb\endcsname{\color{black}}%
      \expandafter\def\csname LTa\endcsname{\color{black}}%
      \expandafter\def\csname LT0\endcsname{\color[rgb]{1,0,0}}%
      \expandafter\def\csname LT1\endcsname{\color[rgb]{0,1,0}}%
      \expandafter\def\csname LT2\endcsname{\color[rgb]{0,0,1}}%
      \expandafter\def\csname LT3\endcsname{\color[rgb]{1,0,1}}%
      \expandafter\def\csname LT4\endcsname{\color[rgb]{0,1,1}}%
      \expandafter\def\csname LT5\endcsname{\color[rgb]{1,1,0}}%
      \expandafter\def\csname LT6\endcsname{\color[rgb]{0,0,0}}%
      \expandafter\def\csname LT7\endcsname{\color[rgb]{1,0.3,0}}%
      \expandafter\def\csname LT8\endcsname{\color[rgb]{0.5,0.5,0.5}}%
    \else
      \def\colorrgb#1{\color{black}}%
      \def\colorgray#1{\color[gray]{#1}}%
      \expandafter\def\csname LTw\endcsname{\color{white}}%
      \expandafter\def\csname LTb\endcsname{\color{black}}%
      \expandafter\def\csname LTa\endcsname{\color{black}}%
      \expandafter\def\csname LT0\endcsname{\color{black}}%
      \expandafter\def\csname LT1\endcsname{\color{black}}%
      \expandafter\def\csname LT2\endcsname{\color{black}}%
      \expandafter\def\csname LT3\endcsname{\color{black}}%
      \expandafter\def\csname LT4\endcsname{\color{black}}%
      \expandafter\def\csname LT5\endcsname{\color{black}}%
      \expandafter\def\csname LT6\endcsname{\color{black}}%
      \expandafter\def\csname LT7\endcsname{\color{black}}%
      \expandafter\def\csname LT8\endcsname{\color{black}}%
    \fi
  \fi
    \setlength{\unitlength}{0.0500bp}%
    \ifx\gptboxheight\undefined%
      \newlength{\gptboxheight}%
      \newlength{\gptboxwidth}%
      \newsavebox{\gptboxtext}%
    \fi%
    \setlength{\fboxrule}{0.5pt}%
    \setlength{\fboxsep}{1pt}%
    \definecolor{tbcol}{rgb}{1,1,1}%
\begin{picture}(4320.00,4320.00)%
    \gplgaddtomacro\gplbacktext{%
      \csname LTb\endcsname
      \put(-132,440){\makebox(0,0)[r]{\strut{}$-2$}}%
      \csname LTb\endcsname
      \put(-132,1355){\makebox(0,0)[r]{\strut{}$-1$}}%
      \csname LTb\endcsname
      \put(-132,2270){\makebox(0,0)[r]{\strut{}$0$}}%
      \csname LTb\endcsname
      \put(-132,3184){\makebox(0,0)[r]{\strut{}$1$}}%
      \csname LTb\endcsname
      \put(-132,4099){\makebox(0,0)[r]{\strut{}$2$}}%
      \csname LTb\endcsname
      \put(270,220){\makebox(0,0){\strut{}$\frac{1}{4}$}}%
      \csname LTb\endcsname
      \put(540,220){\makebox(0,0){\strut{}$\frac{1}{2}$}}%
      \csname LTb\endcsname
      \put(1080,220){\makebox(0,0){\strut{}$1$}}%
      \csname LTb\endcsname
      \put(2160,220){\makebox(0,0){\strut{}$2$}}%
      \csname LTb\endcsname
      \put(4319,220){\makebox(0,0){\strut{}$4$}}%
    }%
    \gplgaddtomacro\gplfronttext{%
    }%
    \gplbacktext
    \put(0,0){\includegraphics[width={216.00bp},height={216.00bp}]{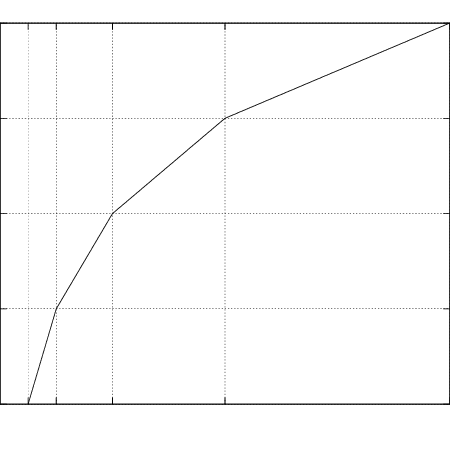}}%
    \gplfronttext
  \end{picture}%
\endgroup

%% file: coarse_approximation.tex
\subsection{Coarse Approximation}

As has been well-established by prior work, given a suitable choice for the constant $c$, the value of $y$ in line \ref{alg:FRSR:line:y} of Algorithm \ref{alg:FRSR} gives a coarse approximation to $1/ \sqrt{x}$. That is, multiplication of $X$ by $-\frac{1}{2}$ in pseudolog space approximates raising $x$ to the power $-\frac{1}{2}$. By the same token, we can approximate other fractional powers of $x$ by multiplying $X$ by the appropriate fraction.\footnote{On many platforms this fractional multiply can be achieved using a small number of low-cost integer operations - for example, when using \texttt{gcc} to compile for a 64-bit x86 target, the integer expression \texttt{X/3} typically generates an integer multiply followed by a shift operation.}

In the present work, for reasons which will become clear, only rational powers of $x$ will be considered, and so the target function to be approximated can be written as
\begin{equation} \label{eq:requirements}
  f(x)=x^{-\frac{a}{b}}\,,  \quad  x \in \mathbb{R^+}\,, \enskip a,b \in \mathbb{Z}^+\,, \enskip gcd(a,b)=1\,.
\end{equation}
Note that we are assuming the fraction $\frac{a}{b}$ has been reduced to its lowest terms, and that only negative powers of $x$ are being considered. (Positive rational powers can then be obtained by multiplying by a suitable integer power of $x$, which can be achieved using only multiply instructions.) Note also that we are considering the domain of $f$ to include all positive real numbers, and that all calculations will be considered exact. (This assumption obviously does not hold in a practical implementation, where only finite precision is available, and where the optimal choices for the values of the constants may deviate from the analytically derived ones. We will return to this problem in section \ref{sec:implementation}.)

To generalise the coarse approximation in Algorithm \ref{alg:FRSR}, we replace the fractions $\frac{c}{2}$ and $\frac{1}{2}$ in line \ref{alg:FRSR:line:Y} with $\frac{c}{b}$ and $\frac{a}{b}$, respectively, so that the following linear relationship holds between $X$ and $Y$:
\begin{equation} \label{eq:line-XY}
  aX+bY=c\,.
\end{equation}
In justification of this, since $X \approx \log_2{x}$ and $Y \approx \log_2{y}$, then
\[  y \approx 2^Y = 2^{\frac{c-aX}{b}} \approx 2^{\frac{c}{b}} 2^{-\frac{a}{b} \log_2 x} = 2^{\frac{c}{b}} x^{-\frac{a}{b}}\,,  \]
so that our coarse approximation $y$ is roughly the target function $x^{-\frac{a}{b}}$ scaled by the constant value $2^{\frac{c}{b}}$. Setting $c$ to zero would provide a simple way to have $y \approx x^{-\frac{a}{b}}$, though as the FRSR case suggests, this is not generally optimal. We show later how to choose an optimal value for $c$.

%% file: refined_approximation.tex
\subsection{Refined Approximation} \label{sec:refined-approximation}

We begin by formalising the definition of the pseudolog function $L(x)$. Given $x \in \mathbb{R^+}$, let
\begin{equation} \label{eq:Em(x)}
  E_x = \floor{\log_2 {x}}\,, \quad m_x = 2^{-E_x} x\ - 1,
\end{equation}
with $\floor{\ }$ denoting the floor function, so that $E_x \in \mathbb{Z}$ and $0\leqslant m_x < 1$. We will refer to $E_x$ and $m_x$ as the \textit{exponent} and \textit{mantissa}, respectively, of $x$, consistent with their counterparts in the terminology of floating point numbers. We may now define $L(x)$ using
\begin{equation} \label{eq:L(x)}
  L(x)=E_x+m_x\,.
\end{equation}
The inverse function $L^{-1}(X)$ can be computed as follows. Given $X \in \mathbb{R}$, let
\begin{equation} \label{eq:Em(X)}
  E_{x} = \floor{X}\,, \quad m_{x} = X - E_{x}\,.
\end{equation}
Then
\begin{equation} \label{eq:L_inv(X)}
  L^{-1}(X) = 2^{E_{x}} (1+m_{x})\,.
\end{equation}
To develop a refined approximation, first let $g(x)$ be the factor by which we would need to scale the coarse approximation $y$ to exactly match the target function $f$:
\[  g(x) = \frac{f(x)}{y(x)} = x^{-\frac{a}{b}}y^{-\frac{b}{b}} = z^{-\frac{1}{b}}\,,  \]
where the auxiliary function $z(x)$ is defined by
\begin{equation} \label{eq:z}
  z(x) = x^a (y(x))^b\,.
\end{equation}
Intuitively, $z$ measures the goodness of fit of $y$ to the target function $f$, with a value of $1$ indicating a perfect match, leading to the notion that the deviation of $z$ from $1$ can be used in computing a correction factor.

If we had an efficient way to compute $g$ exactly, we could generate the precise value of the target function $f$. This is of course an unreasonable expectation, and in practice the best we can do is to approximate $g(x)$ via some method. A polynomial in $x$ might seem an obvious choice, but we run into a problem: the domain is infinite. However, it can easily be shown that $z$ is bounded. Indeed, by decomposing $x$ and $y$ into their respective exponents and mantissas using transformations \eqref{eq:Em(x)} - \eqref{eq:L_inv(X)}, and substituting into equation \eqref{eq:z}, we have
\begin{equation} \label{eq:z(X,Y)}
  z = 2^{a E_x + b E_y} (1+m_x)^a (1+m_y)^b \,,
\end{equation}
and since $X-1<E_x\leqslant X$ and $1\leqslant1+m_x<2$, and similarly for $E_y$ and $m_y$, we can apply the linear relationship \eqref{eq:line-XY} to show that
\[  2^{c-a-b} = 2^{a(X-1)+b(Y-1)} 1^a 1^b < z < 2^{aX+bY} 2^a 2^b = 2^{a+b+c}\,,  \]
giving lower and upper bounds for $z$. In the following sections we develop sharp bounds, but for now it suffices to observe that the bounded nature of $z$ will allow us to formulate an approximation for $g$ in terms of $z$. We will use an $n^\text{th}$-degree polynomial $p(z)$ to approximate $g$; our refined approximation will then be
\begin{equation} \label{eq:y-tilde}
  \widetilde{y}(x) = y(x) p(z)\,.
\end{equation}
Notice that this form matches that of the standard FRSR algorithm, where $p$ has degree $1$.

This leads us to the generalised form of the FRSR algorithm, shown in Algorithm \ref{alg:FRGR}, and referred to here as the Fast Reciprocal General Root (FRGR) algorithm.

\begin{algorithm}
\caption{FRGR algorithm}
\label{alg:FRGR}
\begin{algorithmic}[1]
\STATE \textbf{function} FRGR($x,a,b,n$)
\STATE $X = L(x)$                                                                       \label{alg:FRGR:line:X}
\STATE $Y = \frac{c}{b}-\frac{a}{b}X$                                      \label{alg:FRGR:line:Y}
\STATE $y = L^{-1}(Y)$                                                              \label{alg:FRGR:line:y}
\STATE $z = x^ay^b$                                                                 \label{alg:FRGR:line:z}
\STATE \textbf{return} $ yp(z) \quad // \enskip p $ has degree $n$        \label{alg:FRGR:line:return}
\end{algorithmic}
\end{algorithm}

Consider now the relative error $\widetilde{e}$ incurred by the refined approximation $\widetilde{y}$. It is\footnote{For relative error, we use the definition given in many texts on numerical analysis. Many works on the FRSR algorithm use the negative of this definition; it does not affect the analysis in any significant way.}
\begin{equation} \label{eq:e-tilde}
  \widetilde{e}(x) = \frac{f(x)-\widetilde{y}(x)}{f(x)} = 1 - x^{\frac{a}{b}} y p(z) = \frac{z^{-\frac{1}{b}} - p(z)}{z^{-\frac{1}{b}}}\,.
\end{equation}
Note that the rightmost expression here is the relative error incurred when approximating the function $z^{-\frac{1}{b}}$ by the polynomial $p(z)$ on a particular \textit{finite} domain, say $[z_\text{min},z_\text{max}]$. Finding the optimal $p(z)$ is thus a problem which yields to standard minimax theory (see, for example, \citet{fike1968}). The unique solution is given by setting $p(z)$ equal to the appropriate minimax polynomial\footnote{Throughout this paper, unless otherwise stated, we use the term ``minimax polynomial" to mean the polynomial of prescribed maximum degree which minimises the maximum magnitude of the \textit{relative} error, rather than of the absolute error used in some other contexts. }, which can be found using standard numerical methods, such as the Remez Exchange Algorithm \citep{remez1934}.

However, determination of the minimax coefficients requires that  we know the range of values $[z_\text{min}, z_\text{max}]$ over which $z(x)$ can vary for $x \in (0,\infty)$. This range is not yet fixed, because it is dependent on the value $c$. We now turn our attention to investigating the nature of this dependency, i.e. to the problem of finding expressions for $z_\text{min}$ and $z_\text{max}$ as functions of $c$.

%% file: finding_extrema.tex
\subsection{Finding the extrema of $z(x)$}

As has been done in prior research, we shall make use of a graph to help with the analysis. Here, though, we make the observation that the analysis of critical points of $z$ is made dramatically simpler by plotting the coarse approximation function in pseudolog-pseudolog space since, as a consequence of \eqref{eq:line-XY}, it will always be a straight line. That is, instead of plotting $y$ against $x$, we plot $Y=L(y)$ against $X=L(x)$. We shall find that this reduces much of the problem at hand to the classification of certain crossing points.

Consider a point $P(X,Y)$ lying on the line \eqref{eq:line-XY}. The value of $z$ at $P$, which we shall write as $z|_P$ is given by equation \eqref{eq:z(X,Y)}.

Now, clearly the point $P'(X+b,Y-a)$ must also lie on the line. If the exponents corresponding to $X+b$ and $Y-a$ are labelled respectively $E_x'$ and $E_y'$, and the corresponding mantissas $m_x'$ and $m_y'$, by \eqref{eq:Em(X)} we have
\begin{align*}
  E_x'&=\floor{X+b}=E_x+b\,, \quad\quad m_x'=X+b-E_x'=m_x\,, \\
  E_y'&=\floor{Y-a}=E_y-a\,, \quad\quad m_y'=Y-a-E_y'=m_y\,.
\end{align*}
Thus only the exponents change, and the value of $z$ at $P'$ is
\begin{equation} \label{eq:periodicity}
  z|_{P'} = 2^{a (E_x+b) + b (E_y-a)} (1+m_x)^a (1+m_y)^b = z|_P\,.
\end{equation}
We conclude that the function $z$ is periodic on the line \eqref{eq:line-XY}, with period at most $b$ along the $X$-axis and, equivalently, period at most $a$ along the $Y$-axis. We can thus completely characterise the behaviour of $z$ by examining it over a representative interval, say $X\in[0,b]$.

From Algorithm \ref{alg:FRGR} and equation \eqref{eq:L_inv(X)}, we see that $z$ is a continuous function of $X$, and by the extreme value theorem attains a minimum and a maximum on $X \in [0,b]$. Due to the periodic nature of $z$, these also serve as the minimum and maximum of $z$ over all values of $X$.

We can also see that $z$ must have continuous derivative with respect to $X$ except wherever $X$ or $Y$ is an integer, since these are the points where the exponent of $x$ or $y$ changes. Thus, the minimum and maximum of $z$ must occur either where the line crosses a boundary of an integer $(X,Y)$-grid square, or where $z$ is stationary on the line.

In order to locate the points at which $z$ is stationary, consider again the point $P(X,Y)$. Since $E_x=\floor{X}$ and $E_y=\floor{Y}$, we can say that $P$  lies in the square $[E_x,E_x+1)\times[E_y,E_y+1)$. On this region, the exponents $E_x$ and $E_y$ can be treated as constants - only the mantissas $m_x$ and $m_y$ vary. To find the derivative of $z$ with respect to $X$, we make use of \eqref{eq:line-XY} and \eqref{eq:Em(X)} to establish the following:
\[  \frac{dm_x}{dX}=\frac{dm_y}{dY}=1\,, \quad \frac{dY}{dX}=-\frac{a}{b}\,,  \]
and then apply the chain rule to \eqref{eq:z(X,Y)} to obtain
\[  \frac{dz}{dX} =  2^{aE_x+bE_y} a (1+m_x)^{a-1} (1+m_y)^{b-1} (m_y-m_x) \,.  \]
Since $2^{aE_x+bE_y}>0$, $a>0$, and $m_x,m_y \geqslant 0$, $\frac{dz}{dX}$
is zero if and only if $m_x=m_y$. Thus there is exactly one stationary point of $z$ on the grid square in question, and it is located where the line \eqref{eq:line-XY} crosses the diagonal which joins the corners with $(X,Y)$-coordinates $(E_x,E_y)$ and $(E_x+1,E_y+1)$.

We can summarise these findings by stating that for a given value of $c$, the extrema $z_\text{min}(c)$ and $z_\text{max}(c)$ of $z$ must be contained in the union of the following finite sets:
\begin{align} \label{eq:HVD1}
  H &= \{z: \enspace X \in \mathbb{Z}, \enspace 0 \leqslant X < b \} \,, \nonumber \\
  V &= \{z: \enspace Y \in \mathbb{Z}, \enspace 0 \leqslant X < b \} \,, \nonumber \\
  D &= \{z: \enspace X-Y \in \mathbb{Z}, \enspace 0 \leqslant X < b \}\,,
\end{align}
the notation suggesting points whose horizontal, vertical and diagonal coordinates, respectively, are integers.

We have thus characterised a set of candidates for the extrema of $z$, namely the set $H \cup V \cup D$ in \eqref{eq:HVD1}. We proceed to find expressions for the members of $H$. The members of $V$ can be treated analogously, and the members of $D$ by a slight modification of the same method.

In the relationship defined in \eqref{eq:line-XY}, let $s$ and $t$ be the floor and fractional parts, respectively, of $c$,
\begin{equation} \label{eq:s-t}
  s = \floor{c} \,, \quad  t = c - s \,,
\end{equation}
so that we have
\begin{equation} \label{eq:line-XY-st}
  aX+bY=s+t\,, \quad s \in \mathbb{Z}\,, \quad 0 \leqslant t < 1\,.
\end{equation}
Now let $X$ be an integer, as will be the case for a member of the set $H$. Considering that $s-aX$ and $b$ are both integers, let $q_b$ and $r_b$ be respectively the quotient and remainder upon dividing $s-aX$ by $b$, so that
\[  s-aX = q_b b + r_b\,, \quad q_b,r_b\in\mathbb{Z}\,, \quad 0\leqslant r_b < b \,.  \]
By \eqref{eq:line-XY-st}, therefore,
\[  Y = q_b + \frac{r_b+t}{b}\,,  \]
and, since $0 \leqslant r_b+t < b$, the fractional part of $Y$ (i.e. $Y-\floor{Y}$) is $\frac{r_b+t}{b}$.

We are now in a position to evaluate $z$ for the chosen integer $X$-coordinate. From \eqref{eq:Em(X)} and \eqref{eq:L_inv(X)}, since $X \in \mathbb{Z}$ we have $x=2^X$. We have also broken $Y$ into integer and fraction parts, allowing us to write $y=2^{q_b} (1+\frac{r_b+t}{b})$. Substituting into \eqref{eq:z(X,Y)} we obtain
\begin{equation} \label{eq:z-intX}
  z|_{X \in \mathbb{Z}} = (2^X)^a (2^{q_b})^b \left(1+\frac{r_b+t}{b} \right)^b = 2^{s-r_b} \left(1+\frac{r_b+t}{b} \right)^b\,.
\end{equation}
Note that as $X$ ranges over the representative set $\{0,...,b-1\}$, since $a$ and $b$ are coprime, the remainder $r_b$ varies over the full set of residues modulo $b$. However, $q_b$ has disappeared from the expression for $z$, as a consequence of the periodicity observed in \eqref{eq:periodicity}.

With $s$ and $t$ as in \eqref{eq:s-t}, let us now introduce a function $\zeta_{r,k}(c)$ defined by
\begin{equation} \label{eq:zeta}
  \zeta_{r,k}(c) = 2^{s-r} \left( 1+\frac{r+t}{k} \right)^k\,, \quad c \in \mathbb{R}\,, k \in \mathbb{Z^+},  r \in \mathbb{Z}_k \,.
\end{equation}
We are now able to express \eqref{eq:z-intX} more compactly as
\[  z|_{X \in \mathbb{Z}} = \zeta_{r_b,b}(c)\,,  \]
where $r_b$ is some integer in $\{0,...,b-1 \}$. By completely analogous reasoning, the value of $z$ at a point on the line $\eqref{eq:line-XY}$ such that $Y$ is an integer is found to be
\[  z|_{Y \in \mathbb{Z}} = \zeta_{r_a,a}(c)\,,  \]
with $r_a$ some integer in $\{0,...,a-1 \}$.

A slight variant of this reasoning allows us to evaluate $z$ at points where $X-Y$ is an integer. We first rewrite equation \eqref{eq:line-XY-st} as
\[  aW + \gamma Y=s+t \,,  \]
with $W=X-Y$ and $\gamma = a+b$. Now we let $q_{\gamma}$ and $r_{\gamma}$ be respectively the quotient and remainder on dividing $s-aW$ by $\gamma$, so that
\[  s-aW=q_{\gamma} \gamma +r_{\gamma}\,, \quad q_{\gamma},r_{\gamma} \in \mathbb{Z}\,, \quad 0 \leqslant r_{\gamma} < \gamma \,,  \]
from which we derive
\begin{equation} \label{eq:Y-diag}
  Y=q_{\gamma}+\frac{r_{\gamma}+t}{\gamma}\,,
\end{equation}
and since $0 \leqslant r_{\gamma}+t < \gamma$, this expression decomposes $Y$ into integer and fraction parts. We now note that if $X-Y$ is an integer then we must also have
\begin{equation} \label{eq:X-diag}
  X=q'+\frac{r_{\gamma}+t}{\gamma} \,,
\end{equation}
for some $q' \in \mathbb{Z}$. Since $X$ and $Y$ must satisfy \eqref{eq:line-XY-st},
\[  a \left(q'+\frac{r_{\gamma}+t}{\gamma} \right) + b \left(q_{\gamma}+\frac{r_{\gamma}+t}{\gamma} \right) = s+t\,.  \]
Solving for $q'$ we find
\begin{equation} \label{eq:q-bar}
  q'=\frac{s-bq_{\gamma}-r_{\gamma}}{a}\,.
\end{equation}
Having decomposed both $X$ and $Y$ we may proceed to evaluate $z$ at the point $(X,Y)$. From \eqref{eq:Y-diag} and \eqref{eq:X-diag} we have
\[  x = 2^{q'} \left(1+\frac{r_{\gamma}+t}{\gamma} \right)\,, \quad y = 2^{q_{\gamma}} \left(1+\frac{r_{\gamma}+t}{\gamma} \right)\,,  \]
and therefore
\[  z = 2^{aq'+bq_{\gamma}} \left(1+\frac{r_{\gamma}+t}{\gamma} \right)^{a+b}\,.  \]
Now using \eqref{eq:q-bar} and simplifying, we find
\[  z|_{X-Y \in \mathbb{Z}} = \zeta_{r_{\gamma},\gamma}(c)\,,  \]
where $r_{\gamma} \in \{0,...,\gamma-1 \}$.

This allows us to list the values of the members of sets $H$, $V$ and $D$ explicitly:
\begin{align} \label{eq:HVD2}
  H & = \{ \zeta_{r_b,b}(c):0 \leqslant r_b < b \} \,, \nonumber \\
  V & = \{ \zeta_{r_a,a}(c):0 \leqslant r_a < a \} \,, \nonumber \\
  D & = \{ \zeta_{r_{\gamma},\gamma}(c):0 \leqslant r_{\gamma} < \gamma \}\,.
\end{align}
It remains to determine which of the elements in these sets are in fact the extrema of $z$. We break this task into two parts: first, we identify which set contains $z_\text{min}$, and which contains $z_\text{max}$. We then pick out the extrema from among the members of each set so identified.

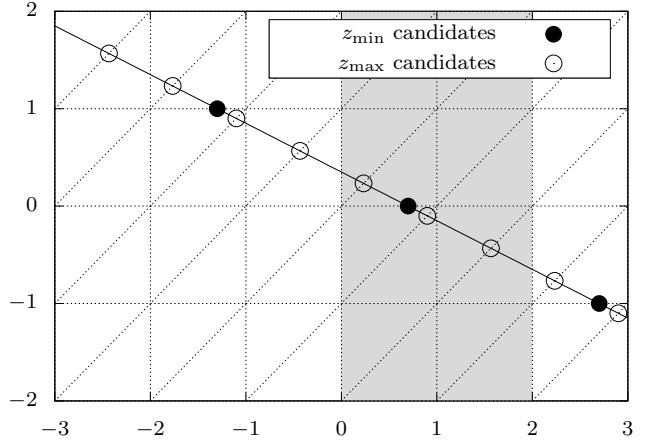
\begin{figure}
  \centering
  \input{candidates}
  \caption{Candidates for extrema of $z$ in FRSR case. The representative interval $X \in [0,2]$ has been shaded.}
  \label{fig:candidates}
\end{figure}

It can be shown (see Appendix \ref{appendix:zeta}) that $\zeta_{r,k}(c)$ is increasing with respect to $k$,
\[  r < k_1 <  k_2 \implies \zeta_{r,k_1}(c) \leqslant \zeta_{r,k_2}(c)\,.  \]
If we now set
\begin{equation} \label{eq:alpha-beta-gamma}
  \alpha = \min(a,b)\,, \quad \beta = \max(a,b)\,, \quad \gamma = a + b \,,
\end{equation}
then we have
\[  0 < \alpha \leqslant \beta < \gamma \,,  \]
and hence for any $r < \alpha$,
\[  \zeta_{r,\alpha}(c) \leqslant \zeta_{r,\beta}(c)\,,  \]
and for any $r < \beta$,
\[  \zeta_{r,\beta}(c) \leqslant \zeta_{r,\gamma}(c)\,.  \]
Comparing with \eqref{eq:HVD2} we conclude that $z_\text{min}$ belongs to $V$ when $a \leqslant b$ and to $H$ when $a \geqslant b$, and that $z_\text{max}$ belongs to the set $D$:
\begin{align}
  z_\text{min}(c) & \in \{ \zeta_{r_{\alpha},\alpha}(c) : 0 \leqslant r_{\alpha} < \alpha \} \,, \label{eq:zmin-candidates} \\
  z_\text{max}(c) & \in \{ \zeta_{r_{\gamma},\gamma}(c) : 0 \leqslant r_{\gamma} < \gamma \} \,. \label{eq:zmax-candidates}
\end{align}
These two sets of candidate points are illustrated for the FRSR case in Figure \ref{fig:candidates} using an arbitrarily chosen value of $0.7$ for $c$. We have shaded the representative interval $[0,b]$, and shown the periodic repeats of the candidate points outside this interval.

The final stage in finding expressions for the extrema is to identify from among the candidates in \eqref{eq:zmin-candidates} which one is in fact $z_\text{min}$, and similarly for \eqref{eq:zmax-candidates} and $z_\text{max}$.

Let $k \geqslant 1$ be a fixed integer, and consider the set of functions
\[  S_k = \left \{ \zeta_{r,k}(c) \right \}_{0 \leqslant r < k}\,.  \]
Clearly, for the simple case $k=1$, since the set contains only the single function $\zeta_{0,1}(c)$, it must serve as both the minimum and maximum member of the set, for all values of $c$. So let us assume that $k \geqslant 2$ in the remainder of this section.

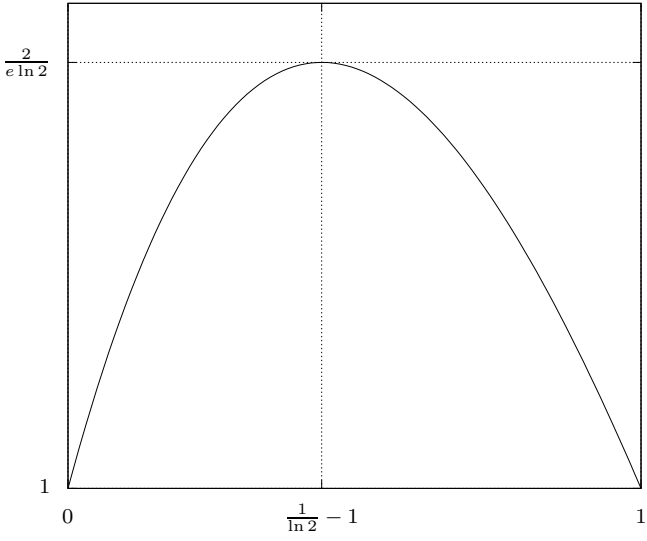
\begin{figure}
  \centering
  \input{eta}
  \caption{The function $\eta(\theta) = 2^{-\theta}(1+\theta)$}
  \label{fig:eta}
\end{figure}

To analyse the general case, consider first the continuous function
\[  \eta(\theta) = 2^{-\theta} \left( 1+\theta \right), \quad \theta \in [0,1]\,,  \]
which is plotted in Figure \ref{fig:eta}.

Elementary calculus enables us to determine that $\eta(\theta)$ is minimum at the endpoints of the interval, is maximum at $\theta = \frac{1}{\ln 2}-1=\bar{\theta}$, say, and that it is strictly increasing for $\theta < \bar{\theta}$ and strictly decreasing for $\theta > \bar{\theta}$. By restricting $\theta$ to the values $\frac{r}{k}$ with $0 \leqslant r < k$ it becomes clear that $\eta(\frac{r}{k})$ is minimum only when $r=0$.

As for the maximum, we note that $\frac{r}{k}$ cannot equal $\bar{\theta}$, the latter being irrational, hence the maximum must occur for one of the two points straddling $\bar{\theta}$ (and there must be two, since $0 < \bar{\theta} < \frac{1}{2} \leqslant \frac{k-1}{k}$). It cannot be both; for if $\eta(\frac{r}{k}) = \eta(\frac{r+1}{k})$, then by expanding and rearranging we have
\[  2^{1/k} = \frac{r+k+1}{r+k}\,,  \]
and clearly if $k \geqslant 2$ then the left side is irrational while the right side is rational, which is impossible.

Hence we can say that there exists a unique $\bar{r} \in \mathbb{Z}_k$ for which $\eta(\frac{\bar{r}}{k})$ is maximal, whence
\begin{align}
  0 \leqslant r < r' \leqslant \bar{r} & \implies \eta \left( \frac{r}{k} \right) < \eta \left( \frac{r'}{k} \right) \,, \label{eq:eta-ineq0} \\
  \bar{r} \leqslant r < r' < k & \implies \eta \left( \frac{r}{k} \right) > \eta \left( \frac{r'}{k} \right) \,.   \label{eq:eta-ineq1}
\end{align}
Now define a new function $\widehat{\zeta}_{r,k}(t)$ closely related to $\zeta_{r,k}(c)$ as follows:
\[  \widehat{\zeta}_{r,k}(t) = 2^{-r/k} \left( 1+\frac{r+t}{k} \right), \quad 0 \leqslant t \leqslant 1, \quad 0 \leqslant r < k\,, \]
and define a corresponding set of functions
\begin{equation} \label{eq:set-S-hat}
  \widehat{S}_k = \left \{ \widehat{\zeta}_{r,k}(t) \right \}_{0 \leqslant r < k}\,.
\end{equation}
It is easily verified that
\begin{equation} \label{eq:zeta_hat(0)}
  \widehat{\zeta}_{r,k}(0) = \eta \left( \frac{r}{k} \right)\,,
\end{equation}
and that
\begin{equation} \label{eq:zeta_hat(1)}
  \widehat{\zeta}_{r,k}(1) = 2^{\frac{1}{k}} \eta \left( \frac{r+1}{k} \right)\,.
\end{equation}
Putting $r'=r+1$ into \eqref{eq:eta-ineq0} and using \eqref{eq:zeta_hat(0)}, we find
\begin{align}
  0 \leqslant r < \bar{r} & \implies \widehat{\zeta}_{r+1,k}(0) - \widehat{\zeta}_{r,k}(0) > 0 \,, \nonumber \\
  \bar{r} \leqslant r < k-1 & \implies \widehat{\zeta}_{r+1,k}(0) - \widehat{\zeta}_{r,k}(0) < 0 \,.  \label{eq:zeta_hat(0)_monotonic}
\end{align}
By similar reasoning using \eqref{eq:eta-ineq1} and \eqref{eq:zeta_hat(1)},
\begin{align}
  0 \leqslant r < \bar{r}-1 & \implies \widehat{\zeta}_{r+1,k}(1) - \widehat{\zeta}_{r,k}(1) > 0 \,, \nonumber \\
  \bar{r}-1 \leqslant r < k-1 & \implies \widehat{\zeta}_{r+1,k}(1) - \widehat{\zeta}_{r,k}(1) < 0 \,.  \label{eq:zeta_hat(1)_monotonic}
\end{align}
Noting that $\widehat{\zeta}$ is linear in $t$, i.e.,
\[  \widehat{\zeta}_{r,k}(t) = (1-t) \widehat{\zeta}_{r,k}(0) + t \widehat{\zeta}_{r,k}(1)\,,  \]
for $t \in [0,1]$ we can say that $\widehat{\zeta}_{r,k}(t)$ is a convex combination of $\widehat{\zeta}_{r,k}(0)$ and $\widehat{\zeta}_{r,k}(1)$, so that \eqref{eq:zeta_hat(0)_monotonic} and \eqref{eq:zeta_hat(1)_monotonic} combine to yield
\begin{align}
  0 \leqslant r < \bar{r}-1 & \implies \widehat{\zeta}_{r+1,k}(t) > \widehat{\zeta}_{r,k}(t) \,, \nonumber \\
  \bar{r} \leqslant r < k-1 & \implies \widehat{\zeta}_{r+1,k}(t) < \widehat{\zeta}_{r,k}(t) \,, \label{eq:zeta-hat-partition}
\end{align}
and this partitions the set $\widehat{S}_k$ of \eqref{eq:set-S-hat} into two subsets, each of which contains non-intersecting functions.

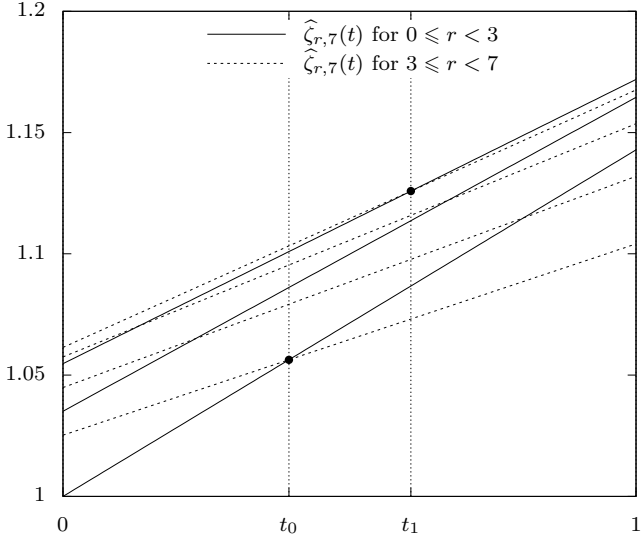
\begin{figure}
  \centering
  \input{zeta-hat}
  \caption{The functions $\widehat{\zeta}_{r,7}(t)$ plotted for $t \in [0,1]$}
  \label{fig:zeta_hat}
\end{figure}

At this point, an example will serve to clarify. Let us consider the case $k=7$, and plot the set of functions $\widehat{S}_7$ for $t \in [0,1]$. The result is shown in Figure \ref{fig:zeta_hat}; the meaning of the symbols $t_0$ and $t_1$ will be explained presently. In this case we have $\bar{r}=3$, and as hinted at by the two different line styles, $\widehat{S}_7$ can be partitioned into two subsets,
\[  \left \{ \widehat{\zeta}_{r,7}(t) \right \}_{0 \leqslant r < 3 } \cup \left \{ \widehat{\zeta}_{r,7}(t) \right \}_{3 \leqslant r < 7 }  \]
each of which contains disjoint line segments and is thus totally ordered.

In the general case, \eqref{eq:zeta-hat-partition} partitions $\widehat{S}_k$ into the subsets
\[  \left \{ \widehat{\zeta}_{r,k}(t) \right \}_{0 \leqslant r < \bar{r} } \cup \left \{ \widehat{\zeta}_{r,k}(t) \right \}_{\bar{r} \leqslant r < k } \,.  \]
Using the total order which \eqref{eq:zeta-hat-partition} induces on each subset, the minimum and maximum members of the first subset are $\widehat{\zeta}_{0,k}(t)$ and $\widehat{\zeta}_{\bar{r}-1,k}(t)$, respectively, and those of the second subset are  $\widehat{\zeta}_{k-1,k}(t)$ and $\widehat{\zeta}_{\bar{r},k}(t)$, respectively. 

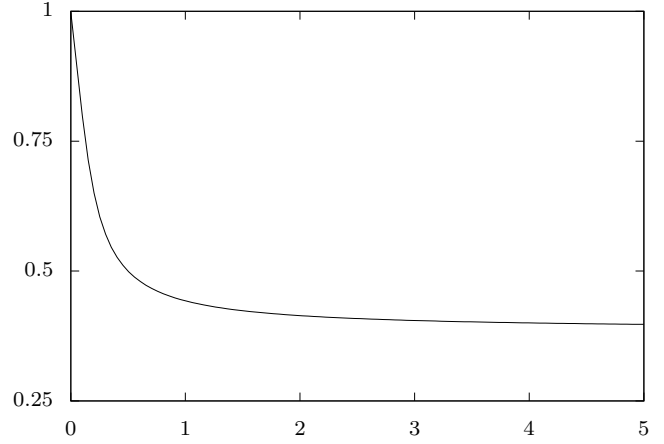
\begin{figure}
  \centering
  \input{t0}
  \caption{The function $t_0(k)$}
  \label{fig:t0}
\end{figure}

The minimum for the full set $\widehat{S}_k$ is obtained by choosing appropriately, depending on the value of $t$, between the pair of subset minima $\widehat{\zeta}_{0,k}(t)$ and $\widehat{\zeta}_{k-1,k}(t)$. Suppose these two functions intersect at $t=t_0$. Setting the two functions equal and solving for $t_0$ yields
\begin{equation} \label{eq:t0}
  t_0(k) = \frac{k-1}{2^{1-\frac{1}{k}}-1}-k\,.
\end{equation}
Although we have stipulated that $k$ be a positive integer, it is instructive to extend the domain of $t_0$ to include all non-negative real numbers. A graph of $t_0(k)$ is shown in Figure \ref{fig:t0}.

Note that this function is undefined when $k$ takes either of the values $0$ or $1$. However, the limit exists in both cases, and we find that
\[  \lim_{k \to 0} t_0(k) = 1\,, \quad \lim_{k \to 1} t_0(k) = \frac{1}{\ln{2}} -1 \,, \]
the first limit being trivial to show, and the second being easily evaluated using l'H$\hat{\text{o}}$pital's rule. For completeness, then, we define $t_0(0)=1$ and $t_0(1)=\frac{1}{\ln{2}}-1$.

In Appendix \ref{appendix:t0}, we demonstrate bounds on $t_0(k)$ (see equation \eqref{eq:t0-bounds}); for present purposes it suffices to note that they imply $t_0(k)$ lies strictly between $0$ and $1$ under our assumption that $k \geqslant 2$.

From \eqref{eq:zeta_hat(0)}, $\widehat{\zeta}_{r,k}(0)$ is minimised by setting $r=0$, and we deduce that the minimum member of $\widehat{S}_k$ is $\widehat{\zeta}_{0,k}(t)$ for values of $t$ less than $t_0(k)$, and is $\widehat{\zeta}_{k-1,k}(t)$ all greater values of $t$:
\begin{equation} \label{eq:min-zeta-hat}
  \min_{0 \leqslant r < k} \widehat{\zeta}_{r,k}(t) =
  \begin{cases}
    \widehat{\zeta}_{0,k}(t)	& \mbox{if } t < t_0(k) \,, \\ 
    \widehat{\zeta}_{k-1,k}(t)	& \mbox{if } t \geqslant t_0(k)\,.
  \end{cases}
\end{equation}
Similarly, for the maximum, we select between the two subset maxima $\widehat{\zeta}_{\bar{r}-1,k}(t)$ and $\widehat{\zeta}_{\bar{r},k}(t)$. Supposing these functions to intersect at $t=t_1$, we shall show that $t_1$ lies strictly between $0$ and $1$.

Since the unique minimum of $\eta(\frac{r}{k})$ occurs for $r=0$, then the maximum must occur for some other value of $r \in \mathbb{Z}_k$, since $k \geqslant 2$. Hence,
\begin{equation} \label{eq:r_bar-bounds}
  1 \leqslant \bar{r} \leqslant k-1 \,.
\end{equation}
We can therefore substitute $r = \bar{r}-1$ into both \eqref{eq:zeta_hat(0)_monotonic} and \eqref{eq:zeta_hat(1)_monotonic} to obtain
\begin{align*}
  \widehat{\zeta}_{\bar{r}-1,k}(0) - \widehat{\zeta}_{\bar{r},k}(0) < 0 \,, \\
  \widehat{\zeta}_{\bar{r}-1,k}(1) - \widehat{\zeta}_{\bar{r},k}(1) > 0 \,.
\end{align*}
By the intermediate value theorem, therefore, the functions $\widehat{\zeta}_{\bar{r}-1,k}(t)$ and $\widehat{\zeta}_{\bar{r},k}(t)$ are equal for some $t$ in the open interval $(0,1)$. By definition, this value is $t_1$; hence $0 < t_1 < 1$.

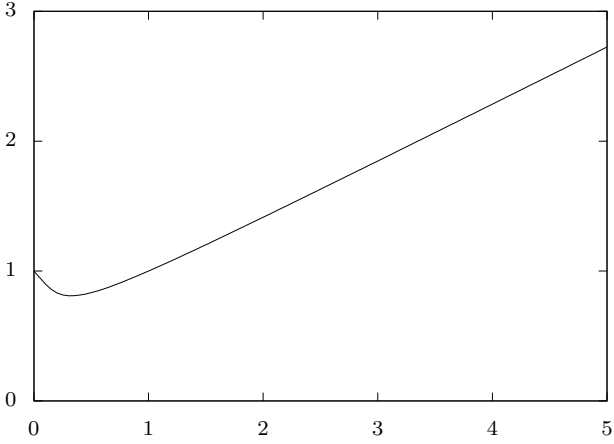
\begin{figure}
  \centering
  \input{phi}
  \caption{The function $\phi(k)$}
  \label{fig:phi}
\end{figure}

By equating the two functions we obtain
\[  2^{-\frac{\bar{r}-1}{k}} \left( 1+\frac{\bar{r}-1+t_1}{k} \right) = 2^{-\frac{\bar{r}}{k}} \left( 1+\frac{\bar{r}+t_1}{k} \right)\,,  \]
which can be rearranged to give
\[  \bar{r}(k) + t_1(k) = \phi(k)\,,  \]
where
\begin{equation} \label{eq:phi}
  \phi(k) = \frac{1}{2^{\frac{1}{k}}-1} - k+1\,.
\end{equation}
A graph of the function $\phi(k)$ is shown in Figure \ref{fig:phi}, where again we have extended the domain to include all non-negative real $k$ (and have defined $\phi(0) = \underset{k \to 0}{\lim} \phi(k) = 1$).

Since $\bar{r} \in \mathbb{Z}$ and $0 < t_1 <1$, clearly $\bar{r}$ and $t_1$ must be the integer and fraction parts, respectively, of $\phi(k)$:
\begin{equation} \label{eq:r_bar-t1}
  \bar{r}(k) = \left\lfloor \phi(k) \right\rfloor\,, \quad  t_1(k) = \phi(k) - \bar{r}(k) \,.
\end{equation}  
Furthermore, since $\widehat{\zeta}_{r,k}(0)$ is maximised by setting $r=\bar{r}$ (by definition of $\bar{r}$), we can deduce that $\widehat{\zeta}_{\bar{r},k}(t)$ is the maximum for values of $t$ less than $t_1$, and that $\widehat{\zeta}_{\bar{r}-1,k}(t)$ is the maximum for all greater values:
\begin{equation} \label{eq:max-zeta-hat}
  \max_{0 \leqslant r < k} \widehat{\zeta}_{r,k}(t) =
  \begin{cases}
    \widehat{\zeta}_{\bar{r},k}(t)	& \mbox{if } t < t_1(k) \,, \\ 
    \widehat{\zeta}_{\bar{r}-1,k}(t)	& \mbox{if } t \geqslant t_1(k)\,.
  \end{cases}
\end{equation}
The intersection points corresponding to $t_0(7)$ and $t_1(7)$ have been highlighted in Figure \ref{fig:zeta_hat}.

At this point, we note that $\zeta$ and $\widehat{\zeta}$ are related via
\[  \zeta_{r,k}(c) = 2^s \left( \widehat{\zeta}_{r,k}(t) \right) ^k \,.  \]
Since $\widehat{\zeta}_{r,k}(t) > 0$, raising it to the positive power $k$ and multiplying by the positive number $2^s$ both preserve monotonicity with respect to $r$ . Hence we can extend results $\eqref{eq:min-zeta-hat}$ and $\eqref{eq:max-zeta-hat}$ from $\widehat{\zeta}$ to $\zeta$ to obtain
\begin{equation} \label{eq:min-zeta}
  \min_{0 \leqslant r < k} \zeta_{r,k}(c) =
  \begin{cases}
    \zeta_{0,k}(c)	& \mbox{if } t < t_0(k) \,, \\ 
    \zeta_{k-1,k}(c) 	& \mbox{if } t \geqslant t_0(k) \,,
  \end{cases}
\end{equation}
\begin{equation} \label{eq:max-zeta}
  \max_{0 \leqslant r < k} \zeta_{r,k}(c) =
  \begin{cases}
    \zeta_{\bar{r},k}(c)	& \mbox{if } t < t_1(k) \,, \\ 
    \zeta_{\bar{r}-1,k}(c)	& \mbox{if } t \geqslant t_1(k) \,.
  \end{cases}
\end{equation}
As the final step, we combine \eqref{eq:min-zeta} with \eqref{eq:zmin-candidates}, and \eqref{eq:max-zeta} with \eqref{eq:zmax-candidates}, to find that
\begin{align}
  z_\text{min}(c) &= \zeta_{r_{\alpha},\alpha}(c) \,,      \label{eq:z_min}  \\
  z_\text{max}(c) &= \zeta_{r_{\gamma},\gamma}(c)\,, \label{eq:z_max}
\end{align}
where the indices $r_{\alpha}$ and $r_{\gamma}$ are given by the piecewise-constant expressions
\begin{equation} \label{eq:r_alpha}
  r_{\alpha} =
  \begin{cases}
    0		& \mbox{if } t < t_0(\alpha) \,, \\
    \alpha-1	& \mbox{if } t \geqslant t_0(\alpha) \,,
  \end{cases}
\end{equation}
\begin{equation} \label{eq:r_gamma}
  r_{\gamma} =
  \begin{cases}
    \bar{r}(\gamma)		& \mbox{if } t < t_1(\gamma) \,, \\
    \bar{r}(\gamma)-1	& \mbox{if } t \geqslant t_1(\gamma)\,.
  \end{cases}
\end{equation}
We now have the desired expressions for $z_\text{min}$ and $z_\text{max}$ in terms of $c$.

For notational compactness, we will henceforth use $t_0, t_1$ and $\bar{r}$ to mean exclusively $t_0(\alpha), t_1(\gamma)$ and $\bar{r}(\gamma)$, respectively.

%% file: candidates.tex
\begingroup
  \inputencoding{cp1252}%
  \fontfamily{Times-New-Roman}%
  \selectfont
  \makeatletter
  \providecommand\color[2][]{%
    \GenericError{(gnuplot) \space\space\space\@spaces}{%
      Package color not loaded in conjunction with
      terminal option `colourtext'%
    }{See the gnuplot documentation for explanation.%
    }{Either use 'blacktext' in gnuplot or load the package
      color.sty in LaTeX.}%
    \renewcommand\color[2][]{}%
  }%
  \providecommand\includegraphics[2][]{%
    \GenericError{(gnuplot) \space\space\space\@spaces}{%
      Package graphicx or graphics not loaded%
    }{See the gnuplot documentation for explanation.%
    }{The gnuplot epslatex terminal needs graphicx.sty or graphics.sty.}%
    \renewcommand\includegraphics[2][]{}%
  }%
  \providecommand\rotatebox[2]{#2}%
  \@ifundefined{ifGPcolor}{%
    \newif\ifGPcolor
    \GPcolorfalse
  }{}%
  \@ifundefined{ifGPblacktext}{%
    \newif\ifGPblacktext
    \GPblacktexttrue
  }{}%
  \let\gplgaddtomacro\g@addto@macro
  \gdef\gplbacktext{}%
  \gdef\gplfronttext{}%
  \makeatother
  \ifGPblacktext
    \def\colorrgb#1{}%
    \def\colorgray#1{}%
  \else
    \ifGPcolor
      \def\colorrgb#1{\color[rgb]{#1}}%
      \def\colorgray#1{\color[gray]{#1}}%
      \expandafter\def\csname LTw\endcsname{\color{white}}%
      \expandafter\def\csname LTb\endcsname{\color{black}}%
      \expandafter\def\csname LTa\endcsname{\color{black}}%
      \expandafter\def\csname LT0\endcsname{\color[rgb]{1,0,0}}%
      \expandafter\def\csname LT1\endcsname{\color[rgb]{0,1,0}}%
      \expandafter\def\csname LT2\endcsname{\color[rgb]{0,0,1}}%
      \expandafter\def\csname LT3\endcsname{\color[rgb]{1,0,1}}%
      \expandafter\def\csname LT4\endcsname{\color[rgb]{0,1,1}}%
      \expandafter\def\csname LT5\endcsname{\color[rgb]{1,1,0}}%
      \expandafter\def\csname LT6\endcsname{\color[rgb]{0,0,0}}%
      \expandafter\def\csname LT7\endcsname{\color[rgb]{1,0.3,0}}%
      \expandafter\def\csname LT8\endcsname{\color[rgb]{0.5,0.5,0.5}}%
    \else
      \def\colorrgb#1{\color{black}}%
      \def\colorgray#1{\color[gray]{#1}}%
      \expandafter\def\csname LTw\endcsname{\color{white}}%
      \expandafter\def\csname LTb\endcsname{\color{black}}%
      \expandafter\def\csname LTa\endcsname{\color{black}}%
      \expandafter\def\csname LT0\endcsname{\color{black}}%
      \expandafter\def\csname LT1\endcsname{\color{black}}%
      \expandafter\def\csname LT2\endcsname{\color{black}}%
      \expandafter\def\csname LT3\endcsname{\color{black}}%
      \expandafter\def\csname LT4\endcsname{\color{black}}%
      \expandafter\def\csname LT5\endcsname{\color{black}}%
      \expandafter\def\csname LT6\endcsname{\color{black}}%
      \expandafter\def\csname LT7\endcsname{\color{black}}%
      \expandafter\def\csname LT8\endcsname{\color{black}}%
    \fi
  \fi
    \setlength{\unitlength}{0.0500bp}%
    \ifx\gptboxheight\undefined%
      \newlength{\gptboxheight}%
      \newlength{\gptboxwidth}%
      \newsavebox{\gptboxtext}%
    \fi%
    \setlength{\fboxrule}{0.5pt}%
    \setlength{\fboxsep}{1pt}%
    \definecolor{tbcol}{rgb}{1,1,1}%
\begin{picture}(4320.00,3600.00)%
    \gplgaddtomacro\gplbacktext{%
    }%
    \gplgaddtomacro\gplfronttext{%
      \csname LTb\endcsname
      \put(3332,3206){\makebox(0,0)[r]{\strut{}$z_\text{min}$ \textrm{candidates}}}%
      \csname LTb\endcsname
      \put(3332,2986){\makebox(0,0)[r]{\strut{}$z_\text{max}$ \textrm{candidates}}}%
      \csname LTb\endcsname
      \put(-132,440){\makebox(0,0)[r]{\strut{}$-2$}}%
      \csname LTb\endcsname
      \put(-132,1175){\makebox(0,0)[r]{\strut{}$-1$}}%
      \csname LTb\endcsname
      \put(-132,1910){\makebox(0,0)[r]{\strut{}$0$}}%
      \csname LTb\endcsname
      \put(-132,2644){\makebox(0,0)[r]{\strut{}$1$}}%
      \csname LTb\endcsname
      \put(-132,3379){\makebox(0,0)[r]{\strut{}$2$}}%
      \csname LTb\endcsname
      \put(0,220){\makebox(0,0){\strut{}$-3$}}%
      \csname LTb\endcsname
      \put(720,220){\makebox(0,0){\strut{}$-2$}}%
      \csname LTb\endcsname
      \put(1440,220){\makebox(0,0){\strut{}$-1$}}%
      \csname LTb\endcsname
      \put(2160,220){\makebox(0,0){\strut{}$0$}}%
      \csname LTb\endcsname
      \put(2879,220){\makebox(0,0){\strut{}$1$}}%
      \csname LTb\endcsname
      \put(3599,220){\makebox(0,0){\strut{}$2$}}%
      \csname LTb\endcsname
      \put(4319,220){\makebox(0,0){\strut{}$3$}}%
    }%
    \gplbacktext
    \put(0,0){\includegraphics[width={216.00bp},height={180.00bp}]{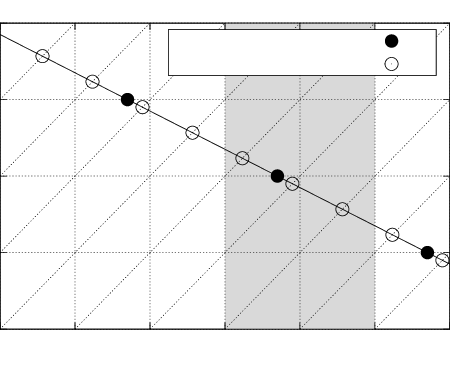}}%
    \gplfronttext
  \end{picture}%
\endgroup

%% file: eta.tex
\begingroup
  \inputencoding{cp1252}%
  \fontfamily{Times-New-Roman}%
  \selectfont
  \makeatletter
  \providecommand\color[2][]{%
    \GenericError{(gnuplot) \space\space\space\@spaces}{%
      Package color not loaded in conjunction with
      terminal option `colourtext'%
    }{See the gnuplot documentation for explanation.%
    }{Either use 'blacktext' in gnuplot or load the package
      color.sty in LaTeX.}%
    \renewcommand\color[2][]{}%
  }%
  \providecommand\includegraphics[2][]{%
    \GenericError{(gnuplot) \space\space\space\@spaces}{%
      Package graphicx or graphics not loaded%
    }{See the gnuplot documentation for explanation.%
    }{The gnuplot epslatex terminal needs graphicx.sty or graphics.sty.}%
    \renewcommand\includegraphics[2][]{}%
  }%
  \providecommand\rotatebox[2]{#2}%
  \@ifundefined{ifGPcolor}{%
    \newif\ifGPcolor
    \GPcolorfalse
  }{}%
  \@ifundefined{ifGPblacktext}{%
    \newif\ifGPblacktext
    \GPblacktexttrue
  }{}%
  \let\gplgaddtomacro\g@addto@macro
  \gdef\gplbacktext{}%
  \gdef\gplfronttext{}%
  \makeatother
  \ifGPblacktext
    \def\colorrgb#1{}%
    \def\colorgray#1{}%
  \else
    \ifGPcolor
      \def\colorrgb#1{\color[rgb]{#1}}%
      \def\colorgray#1{\color[gray]{#1}}%
      \expandafter\def\csname LTw\endcsname{\color{white}}%
      \expandafter\def\csname LTb\endcsname{\color{black}}%
      \expandafter\def\csname LTa\endcsname{\color{black}}%
      \expandafter\def\csname LT0\endcsname{\color[rgb]{1,0,0}}%
      \expandafter\def\csname LT1\endcsname{\color[rgb]{0,1,0}}%
      \expandafter\def\csname LT2\endcsname{\color[rgb]{0,0,1}}%
      \expandafter\def\csname LT3\endcsname{\color[rgb]{1,0,1}}%
      \expandafter\def\csname LT4\endcsname{\color[rgb]{0,1,1}}%
      \expandafter\def\csname LT5\endcsname{\color[rgb]{1,1,0}}%
      \expandafter\def\csname LT6\endcsname{\color[rgb]{0,0,0}}%
      \expandafter\def\csname LT7\endcsname{\color[rgb]{1,0.3,0}}%
      \expandafter\def\csname LT8\endcsname{\color[rgb]{0.5,0.5,0.5}}%
    \else
      \def\colorrgb#1{\color{black}}%
      \def\colorgray#1{\color[gray]{#1}}%
      \expandafter\def\csname LTw\endcsname{\color{white}}%
      \expandafter\def\csname LTb\endcsname{\color{black}}%
      \expandafter\def\csname LTa\endcsname{\color{black}}%
      \expandafter\def\csname LT0\endcsname{\color{black}}%
      \expandafter\def\csname LT1\endcsname{\color{black}}%
      \expandafter\def\csname LT2\endcsname{\color{black}}%
      \expandafter\def\csname LT3\endcsname{\color{black}}%
      \expandafter\def\csname LT4\endcsname{\color{black}}%
      \expandafter\def\csname LT5\endcsname{\color{black}}%
      \expandafter\def\csname LT6\endcsname{\color{black}}%
      \expandafter\def\csname LT7\endcsname{\color{black}}%
      \expandafter\def\csname LT8\endcsname{\color{black}}%
    \fi
  \fi
    \setlength{\unitlength}{0.0500bp}%
    \ifx\gptboxheight\undefined%
      \newlength{\gptboxheight}%
      \newlength{\gptboxwidth}%
      \newsavebox{\gptboxtext}%
    \fi%
    \setlength{\fboxrule}{0.5pt}%
    \setlength{\fboxsep}{1pt}%
    \definecolor{tbcol}{rgb}{1,1,1}%
\begin{picture}(4320.00,4320.00)%
    \gplgaddtomacro\gplbacktext{%
    }%
    \gplgaddtomacro\gplfronttext{%
      \csname LTb\endcsname
      \put(-132,440){\makebox(0,0)[r]{\strut{}$1$}}%
      \csname LTb\endcsname
      \put(-132,3653){\makebox(0,0)[r]{\strut{}$\frac{2}{e \ln{2}}$}}%
      \csname LTb\endcsname
      \put(0,220){\makebox(0,0){\strut{}$0$}}%
      \csname LTb\endcsname
      \put(1912,220){\makebox(0,0){\strut{}$\frac{1}{\ln{2}} - 1$}}%
      \csname LTb\endcsname
      \put(4319,220){\makebox(0,0){\strut{}$1$}}%
    }%
    \gplbacktext
    \put(0,0){\includegraphics[width={216.00bp},height={216.00bp}]{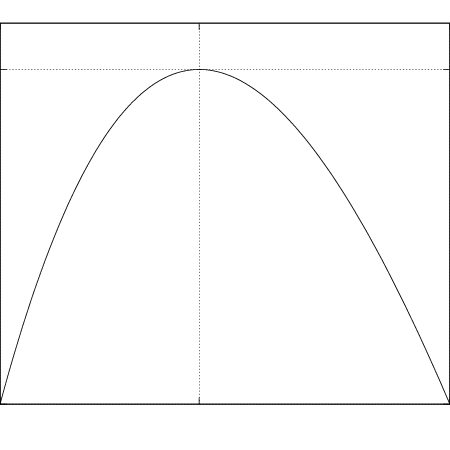}}%
    \gplfronttext
  \end{picture}%
\endgroup

%% file: zeta-hat.tex
\begingroup
  \inputencoding{cp1252}%
  \fontfamily{Times-New-Roman}%
  \selectfont
  \makeatletter
  \providecommand\color[2][]{%
    \GenericError{(gnuplot) \space\space\space\@spaces}{%
      Package color not loaded in conjunction with
      terminal option `colourtext'%
    }{See the gnuplot documentation for explanation.%
    }{Either use 'blacktext' in gnuplot or load the package
      color.sty in LaTeX.}%
    \renewcommand\color[2][]{}%
  }%
  \providecommand\includegraphics[2][]{%
    \GenericError{(gnuplot) \space\space\space\@spaces}{%
      Package graphicx or graphics not loaded%
    }{See the gnuplot documentation for explanation.%
    }{The gnuplot epslatex terminal needs graphicx.sty or graphics.sty.}%
    \renewcommand\includegraphics[2][]{}%
  }%
  \providecommand\rotatebox[2]{#2}%
  \@ifundefined{ifGPcolor}{%
    \newif\ifGPcolor
    \GPcolorfalse
  }{}%
  \@ifundefined{ifGPblacktext}{%
    \newif\ifGPblacktext
    \GPblacktexttrue
  }{}%
  \let\gplgaddtomacro\g@addto@macro
  \gdef\gplbacktext{}%
  \gdef\gplfronttext{}%
  \makeatother
  \ifGPblacktext
    \def\colorrgb#1{}%
    \def\colorgray#1{}%
  \else
    \ifGPcolor
      \def\colorrgb#1{\color[rgb]{#1}}%
      \def\colorgray#1{\color[gray]{#1}}%
      \expandafter\def\csname LTw\endcsname{\color{white}}%
      \expandafter\def\csname LTb\endcsname{\color{black}}%
      \expandafter\def\csname LTa\endcsname{\color{black}}%
      \expandafter\def\csname LT0\endcsname{\color[rgb]{1,0,0}}%
      \expandafter\def\csname LT1\endcsname{\color[rgb]{0,1,0}}%
      \expandafter\def\csname LT2\endcsname{\color[rgb]{0,0,1}}%
      \expandafter\def\csname LT3\endcsname{\color[rgb]{1,0,1}}%
      \expandafter\def\csname LT4\endcsname{\color[rgb]{0,1,1}}%
      \expandafter\def\csname LT5\endcsname{\color[rgb]{1,1,0}}%
      \expandafter\def\csname LT6\endcsname{\color[rgb]{0,0,0}}%
      \expandafter\def\csname LT7\endcsname{\color[rgb]{1,0.3,0}}%
      \expandafter\def\csname LT8\endcsname{\color[rgb]{0.5,0.5,0.5}}%
    \else
      \def\colorrgb#1{\color{black}}%
      \def\colorgray#1{\color[gray]{#1}}%
      \expandafter\def\csname LTw\endcsname{\color{white}}%
      \expandafter\def\csname LTb\endcsname{\color{black}}%
      \expandafter\def\csname LTa\endcsname{\color{black}}%
      \expandafter\def\csname LT0\endcsname{\color{black}}%
      \expandafter\def\csname LT1\endcsname{\color{black}}%
      \expandafter\def\csname LT2\endcsname{\color{black}}%
      \expandafter\def\csname LT3\endcsname{\color{black}}%
      \expandafter\def\csname LT4\endcsname{\color{black}}%
      \expandafter\def\csname LT5\endcsname{\color{black}}%
      \expandafter\def\csname LT6\endcsname{\color{black}}%
      \expandafter\def\csname LT7\endcsname{\color{black}}%
      \expandafter\def\csname LT8\endcsname{\color{black}}%
    \fi
  \fi
    \setlength{\unitlength}{0.0500bp}%
    \ifx\gptboxheight\undefined%
      \newlength{\gptboxheight}%
      \newlength{\gptboxwidth}%
      \newsavebox{\gptboxtext}%
    \fi%
    \setlength{\fboxrule}{0.5pt}%
    \setlength{\fboxsep}{1pt}%
    \definecolor{tbcol}{rgb}{1,1,1}%
\begin{picture}(4320.00,4320.00)%
    \gplgaddtomacro\gplbacktext{%
    }%
    \gplgaddtomacro\gplfronttext{%
      \csname LTb\endcsname
      \put(1811,3926){\makebox(0,0)[l]{\strut{}$\widehat{\zeta}_{r,7}(t)$ for $0 \leqslant r < 3$}}%
      \csname LTb\endcsname
      \put(1811,3706){\makebox(0,0)[l]{\strut{}$\widehat{\zeta}_{r,7}(t)$ for $3 \leqslant r < 7$}}%
      \csname LTb\endcsname
      \put(-132,440){\makebox(0,0)[r]{\strut{}$1$}}%
      \put(-132,1355){\makebox(0,0)[r]{\strut{}$1.05$}}%
      \put(-132,2270){\makebox(0,0)[r]{\strut{}$1.1$}}%
      \put(-132,3184){\makebox(0,0)[r]{\strut{}$1.15$}}%
      \put(-132,4099){\makebox(0,0)[r]{\strut{}$1.2$}}%
      \csname LTb\endcsname
      \put(0,220){\makebox(0,0){\strut{}$0$}}%
      \csname LTb\endcsname
      \put(1703,220){\makebox(0,0){\strut{}$t_0$}}%
      \csname LTb\endcsname
      \put(2622,220){\makebox(0,0){\strut{}$t_1$}}%
      \csname LTb\endcsname
      \put(4319,220){\makebox(0,0){\strut{}$1$}}%
    }%
    \gplbacktext
    \put(0,0){\includegraphics[width={216.00bp},height={216.00bp}]{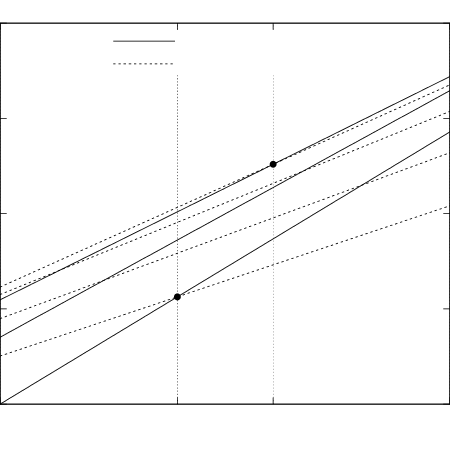}}%
    \gplfronttext
  \end{picture}%
\endgroup

%% file: t0.tex
\begingroup
  \inputencoding{cp1252}%
  \fontfamily{Times-New-Roman}%
  \selectfont
  \makeatletter
  \providecommand\color[2][]{%
    \GenericError{(gnuplot) \space\space\space\@spaces}{%
      Package color not loaded in conjunction with
      terminal option `colourtext'%
    }{See the gnuplot documentation for explanation.%
    }{Either use 'blacktext' in gnuplot or load the package
      color.sty in LaTeX.}%
    \renewcommand\color[2][]{}%
  }%
  \providecommand\includegraphics[2][]{%
    \GenericError{(gnuplot) \space\space\space\@spaces}{%
      Package graphicx or graphics not loaded%
    }{See the gnuplot documentation for explanation.%
    }{The gnuplot epslatex terminal needs graphicx.sty or graphics.sty.}%
    \renewcommand\includegraphics[2][]{}%
  }%
  \providecommand\rotatebox[2]{#2}%
  \@ifundefined{ifGPcolor}{%
    \newif\ifGPcolor
    \GPcolorfalse
  }{}%
  \@ifundefined{ifGPblacktext}{%
    \newif\ifGPblacktext
    \GPblacktexttrue
  }{}%
  \let\gplgaddtomacro\g@addto@macro
  \gdef\gplbacktext{}%
  \gdef\gplfronttext{}%
  \makeatother
  \ifGPblacktext
    \def\colorrgb#1{}%
    \def\colorgray#1{}%
  \else
    \ifGPcolor
      \def\colorrgb#1{\color[rgb]{#1}}%
      \def\colorgray#1{\color[gray]{#1}}%
      \expandafter\def\csname LTw\endcsname{\color{white}}%
      \expandafter\def\csname LTb\endcsname{\color{black}}%
      \expandafter\def\csname LTa\endcsname{\color{black}}%
      \expandafter\def\csname LT0\endcsname{\color[rgb]{1,0,0}}%
      \expandafter\def\csname LT1\endcsname{\color[rgb]{0,1,0}}%
      \expandafter\def\csname LT2\endcsname{\color[rgb]{0,0,1}}%
      \expandafter\def\csname LT3\endcsname{\color[rgb]{1,0,1}}%
      \expandafter\def\csname LT4\endcsname{\color[rgb]{0,1,1}}%
      \expandafter\def\csname LT5\endcsname{\color[rgb]{1,1,0}}%
      \expandafter\def\csname LT6\endcsname{\color[rgb]{0,0,0}}%
      \expandafter\def\csname LT7\endcsname{\color[rgb]{1,0.3,0}}%
      \expandafter\def\csname LT8\endcsname{\color[rgb]{0.5,0.5,0.5}}%
    \else
      \def\colorrgb#1{\color{black}}%
      \def\colorgray#1{\color[gray]{#1}}%
      \expandafter\def\csname LTw\endcsname{\color{white}}%
      \expandafter\def\csname LTb\endcsname{\color{black}}%
      \expandafter\def\csname LTa\endcsname{\color{black}}%
      \expandafter\def\csname LT0\endcsname{\color{black}}%
      \expandafter\def\csname LT1\endcsname{\color{black}}%
      \expandafter\def\csname LT2\endcsname{\color{black}}%
      \expandafter\def\csname LT3\endcsname{\color{black}}%
      \expandafter\def\csname LT4\endcsname{\color{black}}%
      \expandafter\def\csname LT5\endcsname{\color{black}}%
      \expandafter\def\csname LT6\endcsname{\color{black}}%
      \expandafter\def\csname LT7\endcsname{\color{black}}%
      \expandafter\def\csname LT8\endcsname{\color{black}}%
    \fi
  \fi
    \setlength{\unitlength}{0.0500bp}%
    \ifx\gptboxheight\undefined%
      \newlength{\gptboxheight}%
      \newlength{\gptboxwidth}%
      \newsavebox{\gptboxtext}%
    \fi%
    \setlength{\fboxrule}{0.5pt}%
    \setlength{\fboxsep}{1pt}%
    \definecolor{tbcol}{rgb}{1,1,1}%
\begin{picture}(4320.00,3600.00)%
    \gplgaddtomacro\gplbacktext{%
    }%
    \gplgaddtomacro\gplfronttext{%
      \csname LTb\endcsname
      \put(-132,440){\makebox(0,0)[r]{\strut{}$0.25$}}%
      \put(-132,1420){\makebox(0,0)[r]{\strut{}$0.5$}}%
      \put(-132,2399){\makebox(0,0)[r]{\strut{}$0.75$}}%
      \put(-132,3379){\makebox(0,0)[r]{\strut{}$1$}}%
      \put(0,220){\makebox(0,0){\strut{}$0$}}%
      \put(863,220){\makebox(0,0){\strut{}$1$}}%
      \put(1727,220){\makebox(0,0){\strut{}$2$}}%
      \put(2591,220){\makebox(0,0){\strut{}$3$}}%
      \put(3455,220){\makebox(0,0){\strut{}$4$}}%
      \put(4319,220){\makebox(0,0){\strut{}$5$}}%
    }%
    \gplbacktext
    \put(0,0){\includegraphics[width={216.00bp},height={180.00bp}]{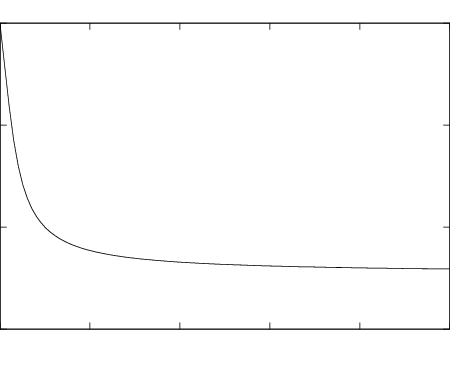}}%
    \gplfronttext
  \end{picture}%
\endgroup

%% file: phi.tex
\begingroup
  \inputencoding{cp1252}%
  \fontfamily{Times-New-Roman}%
  \selectfont
  \makeatletter
  \providecommand\color[2][]{%
    \GenericError{(gnuplot) \space\space\space\@spaces}{%
      Package color not loaded in conjunction with
      terminal option `colourtext'%
    }{See the gnuplot documentation for explanation.%
    }{Either use 'blacktext' in gnuplot or load the package
      color.sty in LaTeX.}%
    \renewcommand\color[2][]{}%
  }%
  \providecommand\includegraphics[2][]{%
    \GenericError{(gnuplot) \space\space\space\@spaces}{%
      Package graphicx or graphics not loaded%
    }{See the gnuplot documentation for explanation.%
    }{The gnuplot epslatex terminal needs graphicx.sty or graphics.sty.}%
    \renewcommand\includegraphics[2][]{}%
  }%
  \providecommand\rotatebox[2]{#2}%
  \@ifundefined{ifGPcolor}{%
    \newif\ifGPcolor
    \GPcolorfalse
  }{}%
  \@ifundefined{ifGPblacktext}{%
    \newif\ifGPblacktext
    \GPblacktexttrue
  }{}%
  \let\gplgaddtomacro\g@addto@macro
  \gdef\gplbacktext{}%
  \gdef\gplfronttext{}%
  \makeatother
  \ifGPblacktext
    \def\colorrgb#1{}%
    \def\colorgray#1{}%
  \else
    \ifGPcolor
      \def\colorrgb#1{\color[rgb]{#1}}%
      \def\colorgray#1{\color[gray]{#1}}%
      \expandafter\def\csname LTw\endcsname{\color{white}}%
      \expandafter\def\csname LTb\endcsname{\color{black}}%
      \expandafter\def\csname LTa\endcsname{\color{black}}%
      \expandafter\def\csname LT0\endcsname{\color[rgb]{1,0,0}}%
      \expandafter\def\csname LT1\endcsname{\color[rgb]{0,1,0}}%
      \expandafter\def\csname LT2\endcsname{\color[rgb]{0,0,1}}%
      \expandafter\def\csname LT3\endcsname{\color[rgb]{1,0,1}}%
      \expandafter\def\csname LT4\endcsname{\color[rgb]{0,1,1}}%
      \expandafter\def\csname LT5\endcsname{\color[rgb]{1,1,0}}%
      \expandafter\def\csname LT6\endcsname{\color[rgb]{0,0,0}}%
      \expandafter\def\csname LT7\endcsname{\color[rgb]{1,0.3,0}}%
      \expandafter\def\csname LT8\endcsname{\color[rgb]{0.5,0.5,0.5}}%
    \else
      \def\colorrgb#1{\color{black}}%
      \def\colorgray#1{\color[gray]{#1}}%
      \expandafter\def\csname LTw\endcsname{\color{white}}%
      \expandafter\def\csname LTb\endcsname{\color{black}}%
      \expandafter\def\csname LTa\endcsname{\color{black}}%
      \expandafter\def\csname LT0\endcsname{\color{black}}%
      \expandafter\def\csname LT1\endcsname{\color{black}}%
      \expandafter\def\csname LT2\endcsname{\color{black}}%
      \expandafter\def\csname LT3\endcsname{\color{black}}%
      \expandafter\def\csname LT4\endcsname{\color{black}}%
      \expandafter\def\csname LT5\endcsname{\color{black}}%
      \expandafter\def\csname LT6\endcsname{\color{black}}%
      \expandafter\def\csname LT7\endcsname{\color{black}}%
      \expandafter\def\csname LT8\endcsname{\color{black}}%
    \fi
  \fi
    \setlength{\unitlength}{0.0500bp}%
    \ifx\gptboxheight\undefined%
      \newlength{\gptboxheight}%
      \newlength{\gptboxwidth}%
      \newsavebox{\gptboxtext}%
    \fi%
    \setlength{\fboxrule}{0.5pt}%
    \setlength{\fboxsep}{1pt}%
    \definecolor{tbcol}{rgb}{1,1,1}%
\begin{picture}(4320.00,3600.00)%
    \gplgaddtomacro\gplbacktext{%
    }%
    \gplgaddtomacro\gplfronttext{%
      \csname LTb\endcsname
      \put(-132,440){\makebox(0,0)[r]{\strut{}$0$}}%
      \put(-132,1420){\makebox(0,0)[r]{\strut{}$1$}}%
      \put(-132,2399){\makebox(0,0)[r]{\strut{}$2$}}%
      \put(-132,3379){\makebox(0,0)[r]{\strut{}$3$}}%
      \put(0,220){\makebox(0,0){\strut{}$0$}}%
      \put(863,220){\makebox(0,0){\strut{}$1$}}%
      \put(1727,220){\makebox(0,0){\strut{}$2$}}%
      \put(2591,220){\makebox(0,0){\strut{}$3$}}%
      \put(3455,220){\makebox(0,0){\strut{}$4$}}%
      \put(4319,220){\makebox(0,0){\strut{}$5$}}%
    }%
    \gplbacktext
    \put(0,0){\includegraphics[width={216.00bp},height={180.00bp}]{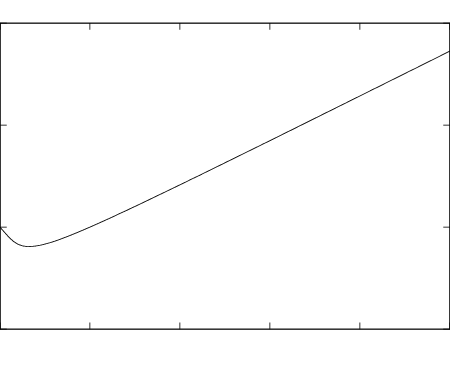}}%
    \gplfronttext
  \end{picture}%
\endgroup

%% file: optimal_c.tex
\subsection{Method for choosing an optimal value for $c$} \label{sec:optimal-c}

With $n$ the degree of the approximating polynomial $p(z)$ in \eqref{eq:y-tilde}, our optimisation problem now requires the determination of $n+2$ values: the $n+1$ coefficients of $p$, and the number $c$. It may appear as though these two problems are inextricably linked. Fortunately, however, there is a method of optimising $c$ independently of $p$.

Given some value of $c$, consider a linear remapping of $z$ given by
\[  \widehat{z}(z)=\frac{z}{z_\text{min}(c)}\,,  \]
and note that as $z$ varies over $[z_\text{min},z_\text{max}]$, so $\widehat{z}$ varies over $[1,\rho]$ with
\begin{equation} \label{eq:rho}
  \rho(c)=\frac{z_\text{max}}{z_\text{min}}\,.
\end{equation}
Note also that the error function in equation \eqref{eq:e-tilde} can be expressed as
\[  \widetilde{e} = \frac{(z_\text{min} \widehat{z})^{-\frac{1}{b}} - p(z_\text{min} \widehat{z})} {(z_\text{min} \widehat{z})^{-\frac{1}{b}}}
                = \frac{\widehat{z}^{-\frac{1}{b}} - q(\widehat{z})} {\widehat{z}^{-\frac{1}{b}}}\,,  \]
where $q(\widehat{z})$ is a polynomial of degree $n$. In other words, the error is the same as that incurred by approximating the function $\widehat{z}^{-\frac{1}{b}}$ by a polynomial of degree $n$ over the interval $[1,\rho]$. Furthermore, the derivative of $\widetilde{e}$ with respect to $\widehat{z}$ is
\[  \frac{d \widetilde{e}}{d \widehat{z}} = -\frac{1}{b} \widehat{z}^{\frac{1}{b}-1} (q(\widehat{z}) + b \widehat{z} q'(\widehat{z}))\,,  \]
which, since $\widehat{z}>0$, is zero only when the polynomial $q(\widehat{z}) + b \widehat{z} q'(\widehat{z})$ is zero. Since the degree of this polynomial is at most $n$, the error function has at most $n$ stationary points.

The Chebyshev Alternation Theorem (see, for example, \citet{fike1968}\footnote{The theorem is usually stated and proved for minimax absolute error approximations. The cited reference also states the theorem in the relative error context suitable for our purposes. The conditions this latter version requires are fulfilled by our scenario.}) tells us that for optimal $q$, the error curve $\widetilde{e}$ must exhibit at least $n+2$ greatest deviations from zero, of equal magnitudes and alternating signs. Since $\widetilde{e}$ has continuous derivative with respect to $\widehat{z}$ on $(0,\infty)$, each point of greatest deviation must either be a stationary point or an endpoint of the interval. And since the stationary points account for at most $n$ of them, the endpoints must constitute the remaining two, and there must be exactly $n+2$ in total.

Now consider the effect of reducing the interval from $[1,\rho]$ to $[1,\rho_1]$, with $0<\rho_1<\rho$, and suppose that a polynomial $q_1(\widehat{z})$ of degree $n$ serves as a minimax approximation on the reduced interval. Let $\epsilon$ and $\epsilon_1$ denote the resulting minimax errors on the original and reduced intervals, respectively.

Since $q(\widehat{z})$ incurs peak error $\epsilon$ on $[1,\rho]$, it incurs no greater an error on a subset of that interval, and hence $\epsilon_1 \leqslant \epsilon$. On the other hand, if $\epsilon_1 = \epsilon$ then both $q$ and $q_1$ are minimax polynomials on $[1,\rho_1]$, and by the uniqueness of the minimax polynomial (see, for example, \citet{fike1968}) we must have $q_1 \equiv q$.

However, the reduced interval excludes the endpoint $\widehat{z}=\rho$ at which one of the original $n+2$ points of greatest deviation lies, leaving at most $n+1$ such points, so that $q$ fails to satisfy the Chebyshev Alternation Theorem on the reduced interval, which is impossible. Hence we must have $\epsilon_1 < \epsilon$.

We conclude that an optimal interval is one which minimises the ratio $\rho = z_\text{max} / z_\text{min}$. Since $z_\text{min}$ and $z_\text{max}$ are dependent on the value of $c$, this minimisation is performed by an appropriate choice of $c$. \footnote{Strictly, we have not addressed the question of the existence of an optimal interval, but our method will demonstrate its existence by construction.}

%% file: minimising_rho.tex
\subsection{Minimising $\rho(c)$}

Using \eqref{eq:z_min} and \eqref{eq:z_max} to expand equation \eqref{eq:rho}, we find
\begin{equation} \label{eq:rho-expanded}
  \rho(c) = \rho(t) = \frac{ 2^{-r_{\gamma}} \left(1+\frac{r_{\gamma}+t}{\gamma} \right) ^ {\gamma} }
			        {  2^{- r_{\alpha}} \left(1+\frac{r_{\alpha}+t}{\alpha} \right) ^ {\alpha} } \,,
\end{equation}
and we wish to minimise this expression on $t \in [0,1)$. Note that $s$, the integer part of $c$, has cancelled from the expression, implying that any minima we find will be repeated for each
value of $s$.

We begin by taking the derivative of \eqref{eq:rho-expanded} with respect to $t$. As a consequence of \eqref{eq:r_alpha} and \eqref{eq:r_gamma} it will have discontinuities at $t=t_0$ and $t=t_1$, but will be continuous in each resulting subinterval of the domain. We find for the derivative\footnote{Strictly, we are here calculating the \textit{right derivative}, as a result of the precise behaviour of $r_{\alpha}$ at $t=t_0$ and of $r_{\gamma}$ at $t=t_1$.}
\begin{align*}
  \frac{d\rho}{dt} & = \left( \frac{ 2^{-r_{\gamma}} \left( 1+\frac{r_{\gamma}+t}{\gamma} \right) ^ {\gamma-1} }
				      {  2^{ -r_{\alpha}} \alpha \gamma \left( 1+\frac{r_{\alpha}+t}{\alpha} \right) ^ {\alpha+1} } \right)
				( \gamma(r_{\alpha}+t) - \alpha(r_{\gamma}+t) )	\\[5pt]
                              & = K(t) \sigma(t) \,,
\end{align*}
where $K(t)$ is the fraction in large parentheses, and $\sigma(t)$ is the piecewise-linear function
\begin{align}	\label{eq:sigma}
  \sigma(t) & = \gamma(r_{\alpha}+t) - \alpha(r_{\gamma}+t)	\nonumber \\
	       & =  \beta t + \gamma r_{\alpha} - \alpha r_{\gamma} \,.
\end{align}
(For both $K$ and $\sigma$ we must bear in mind the fact that $r_{\alpha}$ and $r_{\gamma}$ are themselves piecewise constant functions of $t$.)

Now, $K(t)$ is positive for all $t\in[0,1)$, and so $\rho(t)$ is increasing or decreasing according as the sign of $\sigma(t)$ is positive or negative, respectively. Furthermore, we note that $\sigma(t)$ is the sum of one monotonically increasing function of $t$, namely $\beta t$, plus two weakly increasing functions of $t$, which are $\gamma r_{\alpha}(t)$ and $-\alpha r_{\gamma}(t)$, and is therefore itself monotonically increasing. Hence there is at most one value of $t$ where $\sigma(t)$ and thus $d \rho / dt$ is either zero or jumps at a discontinuity from negative to positive.

Since $t=0$ implies $r_{\alpha}=0,\, r_{\gamma}=\bar{r}$, we have
\[  \sigma(0) = -\alpha \bar{r} < 0 \,,  \]
and, since $t=1$ implies $r_{\alpha}=\alpha-1,\, r_{\gamma}=\bar{r}-1$,
\[  \sigma(1) = \beta + \gamma (\alpha - 1) - \alpha (\bar{r}-1) = \alpha (\gamma - \bar{r}) > 0\,,  \]
the rightmost inequality in both cases following from equation \eqref{eq:r_bar-bounds}.

Furthermore, since we showed that both $t_0$ and $t_1$, the only values of $t$ for which $\sigma(t)$ fails to be continuous, are strictly less than $1$, $\sigma(t)$ must be positive in a small neighbourhood below $1$. 

We conclude that $\sigma(t)$ transitions from negative to positive in the half-open interval $[0,1)$, and since it is monotonically increasing, there must exist exactly one value of $t \in [0,1)$ which minimises $\rho$. This value is therefore the optimum choice for $t$, and we will call it $t^*$. It remains to determine an expression for it.

First, consider the case $\alpha = 1$. From \eqref{eq:r_alpha} we must have $r_{\alpha} = 0$ for all values of $t$, and \eqref{eq:sigma} reduces to
\[
  \sigma(t) = \beta t - r_{\gamma}
                = \begin{cases}
			\beta t - \bar{r}		& \mbox{if } t < t_1 \,, \\
			\beta t - \bar{r} + 1	& \mbox{if } t \geqslant t_1 \,.
		    \end{cases}
\]
Now, if $t_1 \leqslant  \frac{\bar{r}-1} {\beta}$ we have
\[  \sigma \left( \frac{\bar{r}-1} {\beta} \right) = \beta \hspace{1pt} \frac{\bar{r}-1} {\beta} - \bar{r}+1 = 0 \,,  \]
and hence $t^* = \frac{\bar{r}-1} {\beta}$. Similarly, if $t_1 > \frac {\bar{r}} {\beta}$,
\[  \sigma \left( \frac{\bar{r}} {\beta} \right) = \beta \hspace{1pt} \frac{\bar{r}} {\beta} - \bar{r} = 0 \,,  \]
and so $t^* = \frac{\bar{r}} {\beta}$. In the remaining case, i.e. $\frac{\bar{r}-1} {\beta} < t_1 \leqslant  \frac {\bar{r}} {\beta}$,
\[  \sigma(t_1) = \beta t_1 - \bar{r}+1 > 0\,, \quad \lim_{t \to t_1^-} \sigma(t) = \beta t_1 - \bar{r} \leqslant 0 \,,  \]
and in this case we have $t^* = t_1$. We may succinctly state the results from all three cases thus:
\begin{equation} \label{eq:t_star-alpha-1}
  \alpha = 1 \implies t^* = \textrm{clamp} \left( t_1, \frac{\bar{r}-1} {\beta}, \frac{\bar{r}} {\beta} \right) \,.
\end{equation}
We now turn to the case $\alpha \geqslant 2$. Since \eqref{eq:requirements} requires $a$ and $b$ to be coprime,
\begin{equation} \label{eq:alpha-beta-gamma-bounds}
  \alpha \geqslant 2 \implies \beta \geqslant 3 \implies \gamma \geqslant 5 \,.
\end{equation}
In Appendices \ref{appendix:phi} and \ref{appendix:t0}, we demonstrate bounds on $t_0(\alpha)$ and $\phi(\gamma)$ (see equations \eqref{eq:phi-bounds} and \eqref{eq:t0-bounds}) which, to simplify the algebra, we somewhat loosen here to the following:
\begin{equation} \label{eq:t0-bounds-loose}
  \alpha \geqslant 2 \implies 0.3 < t_0(\alpha) < 0.44 \,,
\end{equation}
\begin{equation} \label{eq:phi-bounds-loose}
  \gamma \geqslant 5 \implies \lambda \gamma + 0.5 < \phi(\gamma) < \lambda \gamma + 0.6 \,,
\end{equation}
where
\begin{equation} \label{eq:lambda}
   \lambda = \frac{1}{\ln{2}} - 1\,.
\end{equation}
We may also observe using \eqref{eq:r_bar-t1} and \eqref{eq:r_gamma} that
\begin{equation} \label{eq:r-gamma-ineqs}
  \phi(\gamma)-2 < \bar{r}-1 \leqslant r_{\gamma} \leqslant \bar{r} \leqslant \phi(\gamma)\,.
\end{equation}
Consider setting $t=t_0$ in \eqref{eq:sigma}. We have
\begin{align*}
  \sigma(t_0) & = \beta t_0 + \gamma r_{\alpha}(t_0) - \alpha r_{\gamma}(t_0) \\
                     & \geqslant \beta t_0 + \gamma(\alpha-1) - \alpha \phi(\gamma) & & \text{using \eqref{eq:r_alpha} and \eqref{eq:r-gamma-ineqs}} \\
                     & > 0.3 \beta + \gamma (\alpha - 1) - \alpha ( \lambda \gamma + 0.6 ) & & \text{using \eqref{eq:t0-bounds-loose} and \eqref{eq:phi-bounds-loose}} \\
                     & = (1-\lambda) \alpha \gamma - 0.9 \alpha - 0.7 \gamma \\
                     & > (1-\lambda)(\alpha \gamma - 1.7 \alpha - 1.3 \gamma) & & \text{using \eqref{eq:lambda}} \\
                     & = (1-\lambda)((\alpha - 1.3)(\gamma - 1.7) - 2.21) \\
                     & \geqslant (1-\lambda)((2 - 1.3)(5 - 1.7) - 2.21) & & \text{using \eqref{eq:alpha-beta-gamma-bounds}} \\
                     & = 0.1 (1-\lambda) \\
                     & > 0 \,.
\end{align*}
Now consider the limiting value of $\sigma(t)$ as $t$ approaches $t_0$ from below, $\underset{t \to t_0^-}{\lim} \sigma(t)$, which we will write as $\sigma(t_0^-)$:
\begin{align*}
  \sigma(t_0^-) & = \beta t_0 + \gamma r_{\alpha}(t_0^-) - \alpha r_{\gamma}(t_0^-) \\
                     & < \beta t_0 - \alpha (\phi(\gamma) - 2) & & \text{using \eqref{eq:r_alpha} and \eqref{eq:r-gamma-ineqs}} \\
                     & < 0.44 \beta - \alpha ( \lambda \gamma + 0.5 -2) & & \text{using \eqref{eq:t0-bounds-loose} and \eqref{eq:phi-bounds-loose}} \\
                     & = - \lambda \alpha \gamma + 1.06 \alpha + 0.44 \gamma \\
                     & < -0.44 \alpha \gamma + 1.06 \alpha + 0.44 \gamma & & \text{using \eqref{eq:lambda}} \\
                     & = 1.06 - (\alpha - 1)(0.44 \gamma - 1.06) \\
                     & \leqslant 1.06 - (2 - 1)(2.2 - 1.06) & & \text{using \eqref{eq:alpha-beta-gamma-bounds}} \\
                     & = -0.08 \\
                     & < 0 \,.
\end{align*}
We have thus established that for the case $\alpha \geqslant 2$, $t_0$ is the unique value of $t$ for which $\sigma(t)$ jumps at a discontinuity from negative to positive, leading us to conclude
\begin{equation} \label{eq:t_star-alpha-2}
  \alpha \geqslant 2 \implies t^* = t_0\,.
\end{equation}
Finally, combining \eqref{eq:t_star-alpha-1} with \eqref{eq:t_star-alpha-2}, we obtain
\begin{equation} \label{eq:t-star}
  t^* = \begin{cases}
               \textrm{clamp} \left( t_1, \frac{\bar{r}-1} {\beta}, \frac{\bar{r}} {\beta} \right)    & \mbox{if } \alpha = 1 \,,  \\
               t_0    & \mbox{if } \alpha \geqslant 2 \,.
             \end{cases}
\end{equation}
For any $s \in \mathbb{Z}$, therefore, the value $s+t^*$ with $t^*$ as in \eqref{eq:t-star} provides an optimal choice for the constant c.

\begin{figure}
  \centering
  \input{rho}
  \caption{The function $\rho(c)=\frac{z_\text{max}(c)}{z_\text{min}(c)}$ for the FRSR case. It is periodic outside the interval $[0,1]$.}
  \label{fig:rho}
\end{figure}
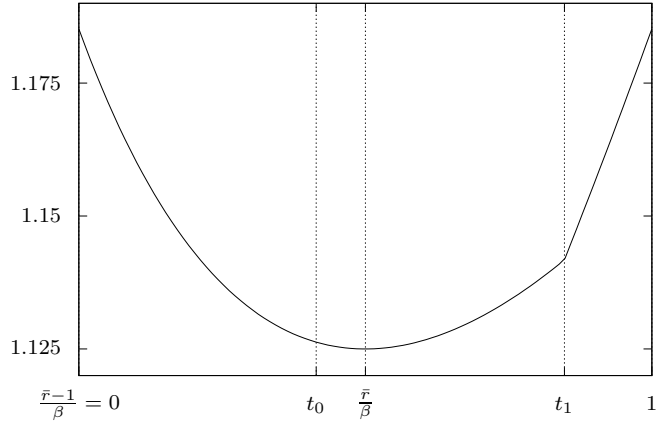

The function $\rho(c)$ corresponding to the FRSR case has been plotted in Figure \ref{fig:rho}. Since $\alpha=1$ for this case, the minimum on $[0,1]$ is located by clamping $t_1$ to the interval $\left[ \frac{\bar{r}-1}{\beta}, \frac{\bar{r}}{\beta} \right]$. The function is periodic outside $[0,1]$, so the minimum repeats for each value of $s$.

%% file: rho.tex
\begingroup
  \inputencoding{cp1252}%
  \fontfamily{Times-New-Roman}%
  \selectfont
  \makeatletter
  \providecommand\color[2][]{%
    \GenericError{(gnuplot) \space\space\space\@spaces}{%
      Package color not loaded in conjunction with
      terminal option `colourtext'%
    }{See the gnuplot documentation for explanation.%
    }{Either use 'blacktext' in gnuplot or load the package
      color.sty in LaTeX.}%
    \renewcommand\color[2][]{}%
  }%
  \providecommand\includegraphics[2][]{%
    \GenericError{(gnuplot) \space\space\space\@spaces}{%
      Package graphicx or graphics not loaded%
    }{See the gnuplot documentation for explanation.%
    }{The gnuplot epslatex terminal needs graphicx.sty or graphics.sty.}%
    \renewcommand\includegraphics[2][]{}%
  }%
  \providecommand\rotatebox[2]{#2}%
  \@ifundefined{ifGPcolor}{%
    \newif\ifGPcolor
    \GPcolorfalse
  }{}%
  \@ifundefined{ifGPblacktext}{%
    \newif\ifGPblacktext
    \GPblacktexttrue
  }{}%
  \let\gplgaddtomacro\g@addto@macro
  \gdef\gplbacktext{}%
  \gdef\gplfronttext{}%
  \makeatother
  \ifGPblacktext
    \def\colorrgb#1{}%
    \def\colorgray#1{}%
  \else
    \ifGPcolor
      \def\colorrgb#1{\color[rgb]{#1}}%
      \def\colorgray#1{\color[gray]{#1}}%
      \expandafter\def\csname LTw\endcsname{\color{white}}%
      \expandafter\def\csname LTb\endcsname{\color{black}}%
      \expandafter\def\csname LTa\endcsname{\color{black}}%
      \expandafter\def\csname LT0\endcsname{\color[rgb]{1,0,0}}%
      \expandafter\def\csname LT1\endcsname{\color[rgb]{0,1,0}}%
      \expandafter\def\csname LT2\endcsname{\color[rgb]{0,0,1}}%
      \expandafter\def\csname LT3\endcsname{\color[rgb]{1,0,1}}%
      \expandafter\def\csname LT4\endcsname{\color[rgb]{0,1,1}}%
      \expandafter\def\csname LT5\endcsname{\color[rgb]{1,1,0}}%
      \expandafter\def\csname LT6\endcsname{\color[rgb]{0,0,0}}%
      \expandafter\def\csname LT7\endcsname{\color[rgb]{1,0.3,0}}%
      \expandafter\def\csname LT8\endcsname{\color[rgb]{0.5,0.5,0.5}}%
    \else
      \def\colorrgb#1{\color{black}}%
      \def\colorgray#1{\color[gray]{#1}}%
      \expandafter\def\csname LTw\endcsname{\color{white}}%
      \expandafter\def\csname LTb\endcsname{\color{black}}%
      \expandafter\def\csname LTa\endcsname{\color{black}}%
      \expandafter\def\csname LT0\endcsname{\color{black}}%
      \expandafter\def\csname LT1\endcsname{\color{black}}%
      \expandafter\def\csname LT2\endcsname{\color{black}}%
      \expandafter\def\csname LT3\endcsname{\color{black}}%
      \expandafter\def\csname LT4\endcsname{\color{black}}%
      \expandafter\def\csname LT5\endcsname{\color{black}}%
      \expandafter\def\csname LT6\endcsname{\color{black}}%
      \expandafter\def\csname LT7\endcsname{\color{black}}%
      \expandafter\def\csname LT8\endcsname{\color{black}}%
    \fi
  \fi
    \setlength{\unitlength}{0.0500bp}%
    \ifx\gptboxheight\undefined%
      \newlength{\gptboxheight}%
      \newlength{\gptboxwidth}%
      \newsavebox{\gptboxtext}%
    \fi%
    \setlength{\fboxrule}{0.5pt}%
    \setlength{\fboxsep}{1pt}%
    \definecolor{tbcol}{rgb}{1,1,1}%
\begin{picture}(4320.00,3456.00)%
    \gplgaddtomacro\gplbacktext{%
    }%
    \gplgaddtomacro\gplfronttext{%
      \csname LTb\endcsname
      \put(-132,640){\makebox(0,0)[r]{\strut{}$1.125$}}%
      \put(-132,1638){\makebox(0,0)[r]{\strut{}$1.15$}}%
      \put(-132,2636){\makebox(0,0)[r]{\strut{}$1.175$}}%
      \csname LTb\endcsname
      \put(0,220){\makebox(0,0){\strut{}$\frac{\bar{r}-1}{\beta}=0$}}%
      \csname LTb\endcsname
      \put(1789,220){\makebox(0,0){\strut{}$t_0$}}%
      \csname LTb\endcsname
      \put(2160,220){\makebox(0,0){\strut{}$\frac{\bar{r}}{\beta}$}}%
      \csname LTb\endcsname
      \put(3660,220){\makebox(0,0){\strut{}$t_1$}}%
      \csname LTb\endcsname
      \put(4319,220){\makebox(0,0){\strut{}$1$}}%
    }%
    \gplbacktext
    \put(0,0){\includegraphics[width={216.00bp},height={172.80bp}]{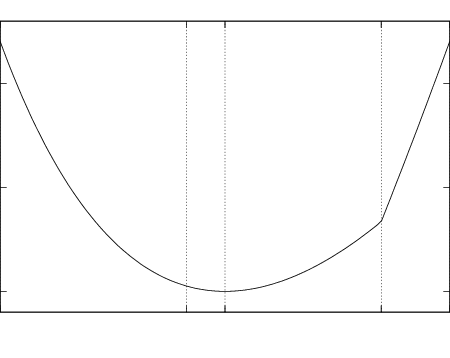}}%
    \gplfronttext
  \end{picture}%
\endgroup

%% file: main_result.tex
\section{Main result}

We have now developed a complete procedure for determining the constants $c, c_0, ..., c_n$ of the FRGR algorithm. We combine equations \eqref{eq:s-t}, \eqref{eq:zeta}, \eqref{eq:alpha-beta-gamma}, \eqref{eq:t0}, \eqref{eq:phi}, \eqref{eq:r_bar-t1}, \eqref{eq:z_min}, \eqref{eq:z_max}, \eqref{eq:r_alpha}, \eqref{eq:r_gamma}  and \eqref{eq:t-star}, together with an assumed numerical procedure \texttt{minimax()} for computing the coefficients of a degree-$n$ minimax polynomial. The result is shown in Algorithm \ref{alg:FRGR-constants}.

\begin{algorithm} [H]
\caption{Determine constants for FRGR}
\label{alg:FRGR-constants}
\begin{algorithmic}[1]
\STATE \textbf{procedure} determine\_FRGR\_constants($a,b,n$)
\STATE $ \alpha = \min(a,b) $
\STATE $ \beta = \max(a,b) $
\STATE $ \gamma = a+b $
\STATE $ t_0 = (\alpha == 1) \enskip ? \enskip \frac{1}{\ln{2}}-1 : \frac{\alpha - 1}{2^{1 - \frac{1}{\alpha}} - 1} - \alpha $
\STATE $ \phi = \frac{1}{2^{\frac{1}{\gamma}}-1} - \gamma + 1 $
\STATE $ \bar{r} = \textrm{floor}(\phi) $
\STATE $ t_1 = \phi - \bar{r} $
\STATE $ t^* = (\alpha == 1) \enskip ? \enskip \textrm{clamp}(t_1, \frac{\bar{r}-1}{\beta}, \frac{\bar{r}}{\beta}) : t_0 $  \label{alg:FRGR-constants:line:t_star}
\STATE $ c = s + t^*, \hspace{1pt} s \in \mathbb{Z} $
\STATE $ r_{\alpha} = (t^* < t_0) \enskip ? \enskip 0 : \alpha - 1 $
\STATE $ r_{\gamma} = (t^* < t_1) \enskip ? \enskip \bar{r} : \bar{r}-1 $
\STATE $ z_\text{min} = 2^{s-r_{\alpha}} \left( 1 + \frac{r_{\alpha}+t^*}{\alpha} \right) ^ {\alpha} $
\STATE $ z_\text{max} = 2^{s-r_{\gamma}} \left( 1 + \frac{r_{\gamma}+t^*}{\gamma} \right) ^ {\gamma} $
\STATE $ p(z) \equiv \textrm{minimax}(z^{-\frac{1}{b}}, n, z_\text{min}, z_\text{max}) $
\STATE \textbf{end procedure}
\end{algorithmic}
\end{algorithm}

%% file: use_of_monics.tex
\section{The use of monic polynomials} \label{sec:monics}

It is common practice to evaluate the polynomial part of the approximation using Horner's rule to reduce the number of floating point operations. For example, in the case $n=3$,
\[  c_3z^3+c_2z^2+c_1z+c_0 = ((c_3z+c_2)z+c_1)z+c_0\,.  \]
In architectures lacking a fused multiply-add (FMA) instruction, this evaluation scheme requires 6 arithmetic operations: 3 multiplies, and 3 adds. In general, a polynomial $p(z)$ of degree $n$ requires $2n$ operations to evaluate in this manner, so that only even numbers of operations arise. It is natural to wonder whether there might be expressions using an odd number of operations which constitute intermediates, in terms of both cost and accuracy. \textit{Monic polynomials}, i.e. polynomials where the leading coefficient is $1$, provide such intermediates.

A monic cubic polynomial, for example, can be evaluated using just 5 operations:
\[  z^3+c_2z^2+c_1z+c_0 =  ((z+c_2)z+c_1)z+c_0\,.  \]
We may also include polynomials with leading coefficient $-1$, which can be evaluated using the same number of operations, by replacing an add with a subtract:
\[  -z^3+c_2z+c_1z+c_0 = ((c_2-z)z+c_1)z+c_0\,.  \]
Suppose in Algorithm \ref{alg:FRGR} that the constant $c$, and hence the approximation interval $[z_\text{min},z_\text{max}]$ has been chosen, and consider the effect of enforcing that the degree-$n$ polynomial $p(z)$ be a signed monic of the form
\[  p(z) = (-z)^n + q(z) \,,  \]
for some polynomial $q$ of degree at most $n-1$, the sign $(-1)^n$ of the leading term of $p$ having been chosen to match that of the $n^\text{th}$-degree general minimax polynomial arising from Algorithm \ref{alg:FRGR-constants}.

If, as in Section \ref{sec:refined-approximation}, $p(z)$ is an approximation for $z^{-\frac{1}{b}}$, the relative error is
\[  \widetilde{e} = \frac{z^{-\frac{1}{b}} - (-z)^n - q(z)}{z^{-\frac{1}{b}}}\,,  \]
the optimisation of which is a \textit{weighted} minimax problem (see, for example, \citet{green2002}): find the polynomial $q$ of degree at most $n-1$ which best approximates $z^{-\frac{1}{b}} - (-z)^n$ using a weight function $z^{\frac{1}{b}}$. Once again, the problem yields to standard minimax theory and the coefficients can be found numerically using a minimax solver.

As before, we can only proceed to find the coefficients once we know the approximation interval. However, in this case, there may not exist a technique that corresponds to the one in Section \ref{sec:optimal-c} which isolates the optimisation of the value $c$ independently from the polynomial coefficients.

On the other hand, it is possible to develop iterative numerical software that converges on the best constant $c$ and simultaneously the best coefficients $c_i$. Empirically, it has been found that for many small values of the input variables $a, b, n$, the optimal value of $c$ is the unique one for which the minimax (general) polynomial happens to be a signed monic. (This is not true in general, a counterexample being $a=b=1, n=3$.)

In the results section, we will use a monic polynomial to demonstrate a reciprocal square root algorithm which is 1 operation faster than the original Quake code while at the same time being roughly twice as accurate. We also exhibit a monic quadratic polynomial for FRSR, showing how, with the addition of a single floating point add instruction, the accuracy can be improved more than 86-fold over the original Quake FRSR.

It may also be observed that, since the monic polynomial of degree $0$ is the constant $1$, a function which simply returns the value of the coarse approximation now fits into our framework, corresponding to the case $p(z) \equiv 1$.

%% file: linear_minimax.tex
\section{Analytic solution of the linear case}

In general, minimax polynomial coefficients do not have closed form expressions, and we must resort to numerical methods to determine them. However, the case of a linear polynomial approximation to $z^{-\frac{1}{b}}$, which occurs in Algorithm \ref{alg:FRGR} for the case $n=1$, is a sufficiently simple one for a closed form solution to exist.

\subsection{\mbox{General linear minimax approximation for $z^{-\frac{1}{b}}$}}

Suppose $c_0 + c_1 z$ is the linear minimax approximation for $z^{-\frac{1}{b}}$ on an interval $[z_\text{min}, z_\text{max}]$. The relative error function has three peaks of equal magnitude and alternating signs, two at the endpoints of the interval and a third, interior point where the error function is stationary. The error can be expressed as
\[  e(z) = 1 - z^{\frac{1}{b}} (c_0 + c_1 z)\,,  \]
with derivative
\[  \frac{de}{dz} = -\frac{1}{b} z^{\frac{1}{b}-1} (c_0 + (b+1) c_1 z)\,.  \]
This can only be zero when $z$ takes the value $z_\text{mid}$, where
\[  z_\text{mid} = -\frac{c_0}{(b+1) c_1}\,,  \]
with equioscillation implying
\begin{equation} \label{eq:linear-equioscillation}
  e(z_\text{min}) = -e(z_\text{mid}) = e(z_\text{max}) = \epsilon\,,
\end{equation}
$\epsilon$ being the minimax error. Solving \eqref{eq:linear-equioscillation} for $c_0, c_1$ and $\epsilon$ leads to
\begin{equation} \label{eq:linear-minimax-solution}
  c_0 = \frac{2T}{U+V}\,, \quad c_1 = \frac{-2}{U+V}\,, \quad \epsilon = \frac{U-V}{U+V}\,,
\end{equation}
with the quantities $T, U, V$ defined as
\begin{align*}
  T & = \frac {z_\text{max}{}^{1+\frac{1}{b}} - z_\text{min}{}^{1+\frac{1}{b}}} {z_\text{max}{}^{\frac{1}{b}} - z_\text{min}{}^{\frac{1}{b}}} \,, \\[5pt]
  U & = b \left( \frac{T}{b+1} \right) ^ {1+\frac{1}{b}} \,, \\[5pt]
  V & = \frac {(z_\text{min} z_\text{max})^{\frac{1}{b}} (z_\text{max} - z_\text{min})} {z_\text{max}{}^{\frac{1}{b}} - z_\text{min}{}^{\frac{1}{b}}} \,.
\end{align*}

\subsection{Application to the FRSR case} \label{sec:app-FRSR}

\citet{walczyk2021} presented an analytic derivation of an optimal set of theoretical constants for the FRSR case. Now that we have developed methods for tackling the more general case, we can use them to quickly derive the FRSR result and confirm their conclusion.

\renewcommand{\thefootnote}{\fnsymbol{footnote}}
The FRSR case corresponds to setting the values $a=1$, $b=2$, $n=1$ in Algorithm \ref{alg:FRGR}, for which the optimal values of $c\footnote[2]{The derivation of the ``magic constant" from the value $c$ will be given in Section \ref{sec:implementation}}, c_0, c_1$ can be obtained using Algorithm \ref{alg:FRGR-constants}. The latter algorithm contains one degree of freedom, the integer $s$; here we choose the value $-1$ so as to correspond with the prior literature. (Choosing any other value would yield an identical value for $\epsilon$, but different values for $c_0$ and $c_1$.)
\renewcommand{\thefootnote}{\arabic{footnote}}

Inserting these values into Algorithm \ref{alg:FRGR-constants}, we obtain
\[  c = -\frac{1}{2}, \quad z_\text{min} = \frac{3}{4}, \quad z_\text{max} = \frac{27}{32} \,,  \]
and the polynomial $p(z) = c_0 + c_1 z$ is then the minimax polynomial of degree $1$ approximating $z^{-\frac{1}{2}}$ on $z \in \left[ \frac{3}{4}, \frac{27}{32} \right]$. Using \eqref{eq:linear-minimax-solution} for the special case $b=2$, we obtain
\begin{align} \label{eq:c0-c1-eps-FRSR}
  c_0 &= \frac { 12(27\sqrt{2}-32) } { \sqrt{10729 - 7242\sqrt{2}} + 9\sqrt{6} } \approx 1.68191391\,, \nonumber \nonumber \\[5pt]
  c_1 &= \frac { 128(4 - 3\sqrt{2}) } { \sqrt{10729 - 7242\sqrt{2}} + 9\sqrt{6} } \approx -0.703952009\,, \nonumber \nonumber \\[5pt]
  \epsilon &=  \frac { \sqrt{10729 - 7242\sqrt{2}} - 9\sqrt{6} } { \sqrt{10729 - 7242\sqrt{2}} + 9\sqrt{6} } \approx 6.50070298 \times 10^{-4} \,.
\end{align}
Our results agree identically with those found by Walczyk et al.

\subsection{Application to the reciprocal function}

It is worth noting in passing that the case of the reciprocal function with a linear minimax polynomial, where $a = b = n = 1$, yields a particularly simple solution under our analysis. For this case, Algorithm \ref{alg:FRGR-constants} together with equation \eqref{eq:linear-minimax-solution} yields the following values (once again choosing the value $-1$ for $s$):
\[  c = \sqrt{2}-2, \quad z_\text{min} = \frac{\sqrt{2}}{2}, \quad z_\text{max} = \frac{1}{8} (3 + 2\sqrt{2}) \,,  \]
\begin{align*}
  c_0 &= \frac {192}{6913} (206\sqrt{2} - 191) \approx 2.78648558\,, \\[3pt]
  c_1 &= \frac {512}{6913} (84\sqrt{2} - 145) \approx -1.94090888 \,, \\[3pt]
  \epsilon &=  \frac {1}{6913} (4481 - 3168\sqrt{2}) \approx 1.11591842 \times 10^{-4}\,.
\end{align*}
We can observe that, even with one fewer multiply instruction than the FRSR algorithm, it can achieve almost $6$ times better accuracy.

%% file: multiple_iterations.tex
\section{Multiple iterations}

The original Quake code \citep{id1998} contained an optional second Newtonian iteration to further refine the result, at the cost of some additional floating point operations. \citet{walczyk2021} used analysis to dramatically improve the accuracy of the 2-iteration version while incurring only a single extra multiply instruction over Quake's 2-iteration version. Here we show how to optimise the coefficients used in each iteration step in the \textit{general} case, further extending the FRGR algorithm. We also show how to remove a multiply instruction from each surplus iteration while retaining the improved accuracy. In the FRSR case, this allows us to keep both the improved accuracy \textit{and} the original execution cost.

\subsection{Generalised multiple iterations} \label{sec:generalised-multiple-iterations}

For this section we will use a slight modification to our notation. We will rewrite the coarse approximation $y$ as $y_0$, and the refined approximation $\widetilde{y}$ as $y_1$; we will then use $y_2$, $y_3$, ... to denote further refined approximations, each having been computed by applying an iteration of the refinement process to the previous approximation. Similarly, we will write $z_0$ in place of $z$ and $p_0$ in place of $p$, and use $z_1$, $p_1$, $z_2$, $p_2$ etc. for their counterparts in subsequent iterations.

By a simple extension of the single-iteration case, a sequence of $m$ refinement steps takes on the following structure:
\begin{align} \label{eq:iterations}
  & z_0 = x^a y_0^b \,, \nonumber \\
  & y_1 = y_0 p_0 (z_0) \,, \nonumber \\
  & z_1 = x^a y_1^b \,, \nonumber \\
  & y_2 = y_1 p_1 (z_1) \,, \nonumber \\
  & ... \nonumber \\
  & z_{m-1} = x^a y_{m-1}^b \,, \nonumber \\
  & y_m = y_{m-1} p_{m-1} (z_{m-1}) \,.
\end{align}
Note that we are free to choose a different degree for each of the $p_i$ if we so wish.

We use a greedy algorithm (a choice which we justify shortly), so that each iteration reduces its own peak relative error to be as small as possible, given the required polynomial degrees. To achieve this, the constant $c$ and refinement polynomial $p_0(z_0)$ are chosen according to Algorithm \ref{alg:FRGR-constants}, and for $i>0$ we choose for $p_i(z_i)$ the minimax polynomial approximation to $z_i^{-\frac{1}{b}}$ on the interval $z_i \in [z_{i;\text{min}}, z_{i;\text{max}}]$. We determine each successive interval as follows.

Given the $i^\text{th}$ approximation $y_i$, we use \eqref{eq:iterations} to calculate $z_i = x^a y_i^b$, and supply this value as the argument to the minimax polynomial $p_i(z_i)$ whose coefficients and domain we assume are already known. Suppose the magnitude of the resulting minimax error is $\epsilon_i$; then for $z_i \in [z_{i;\text{min}},z_{i;\text{max}}]$, we have
\[  -\epsilon_i \leqslant 1 - z_i^{\frac{1}{b}} p_i(z_i) \leqslant \epsilon_i \,,  \]
which we can rearrange as
\[  1-\epsilon_i \leqslant z_i^{\frac{1}{b}} p_i(z_i) \leqslant 1+\epsilon_i \,.  \]
It is easy to show (see Appendix \ref{appendix:eps-bound}) that $1-\epsilon_i > 0$, and hence
\[  (1-\epsilon_i)^b \leqslant z_i p_i(z_i)^b \leqslant (1+\epsilon_i)^b\,.  \]
Now using \eqref{eq:iterations}, we observe that $z_i p_i(z_i)^b = z_{i+1}$, so that this becomes simply
\begin{equation} \label{eq:iter-z-range}
  z_{i+1} \in [(1-\epsilon_i)^b, (1+\epsilon_i)^b]\,,
\end{equation}
which supplies the domain on which we should optimise the polynomial $p_{i+1}(z_{i+1})$.

We justify the use of a greedy algorithm by noting that applying equation \eqref{eq:iter-z-range} to the $(i+1)^\text{st}$ version of the ratio defined in \eqref{eq:rho} yields
\[  \rho_{i+1} = \frac{z_{i+1;\text{max}}}{z_{i+1;\text{min}}} = \left( \frac{1+\epsilon_i}{1-\epsilon_i} \right)^b \,,  \]
and this ratio is clearly minimised by making $\epsilon_i$ as small as possible, which is in turn ensured by choosing for $p_i$ the relevant minimax polynomial. Induction on the number of iterations then completes the argument.

It is worth noting here that if we extend the FRSR algorithm with a second iteration, using a linear polynomial in each iteration as is done in the 2-iteration version of the Quake code, by combining result \eqref{eq:iter-z-range} above with equations \eqref{eq:linear-minimax-solution} from earlier, we find that if $\epsilon_0$ is the analytical minimax error from the first iteration, then that of the second iteration is
\[  \epsilon_1 = \frac{ \left( 1+\frac{\epsilon_0^2}{3} \right)^{\frac{3}{2}} - 1 + \epsilon_0^2 } { \left( 1+\frac{\epsilon_0^2}{3} \right)^{\frac{3}{2}} + 1 - \epsilon_0^2 }\,.  \]
Setting $\epsilon_0$ to the value we obtained in \eqref{eq:c0-c1-eps-FRSR}, we find that the best we can expect from such a 2-iteration FRSR is an error of approximately $3.16943580 \times 10^{-7}$. This value agrees with the one found by \citet{walczyk2021}.

\subsection{Accelerating multiple iterations} \label{sec:acc-multi-iter}

Although Section \ref{sec:generalised-multiple-iterations} specified a recipe for explicitly computing all the coefficients of each polynomial in the sequence of refinements \eqref{eq:iterations}, we actually have some freedom to modify them without changing the final value $y_m$.

To see this, we note that the concatenation of all the steps in \eqref{eq:iterations} yields
\[  y_m = y_0 \prod_{i=0}^{m-1} p_i (z_i)\,,  \]
which contains the product of $m$ polynomials. If we scale each of the $p_i$ by a scaling factor, possibly different for each $i$, then provided the product of the scaling factors is constrained to equal $1$, $y_m$ will be unaltered. Since this scaling scheme has $m$ parameters and one constraint, there are $m-1$ degrees of freedom, which we can use to accelerate the resulting implementation, as we shall see.

However, the modification required to scale a given polynomial is not in general simply a matter of scaling all its coefficients by the same factor, because the value of the polynomial in one such modified iteration affects the argument of the polynomial in the subsequent iteration.

Consider the effect of scaling the coefficients of a single $p_i$ by a non-zero factor $k_i$ to yield a scaled polynomial $p_i'$:
\begin{equation} \label{eq:scale-p_i}
  p_i'(z_i) = k_i p_i(z_i) \,.
\end{equation}
In turn, this scales $y_{i+1}$ to a new value,
\[  y_{i+1}' = y_i k_i p_i(z_i) = k_i y_{i+1}\,,  \]
and hence $z_{i+1}$ to a new value,
\[  z_{i+1}' = x^a y_{i+1}'^b = k_i^b z_{i+1}\,. \]
If we wish the value computed for the polynomial in the next iteration, $p_{i+1}$, to remain unaffected, we must modify its coefficients to counteract the change in its argument, yielding a new polynomial $p_{i+1}'$ where
\begin{equation} \label{eq:scale-p_i+1}
  p_{i+1}' (k_i^b z_{i+1}) = p_{i+1} (z_{i+1})\,.
\end{equation}
Hence the coefficient of the degree-$r$ term in $p_{i+1}'$ should be set equal to the corresponding coefficient in $p_{i+1}$ divided by $k_i^{rb}$.

Together, the transformations \eqref{eq:scale-p_i} and \eqref{eq:scale-p_i+1} give the procedure for scaling the value of a single one of the polynomials by a given constant. We then repeat this procedure for each polynomial we wish to scale. (Note that when scaling the final polynomial $p_{m-1}$ there is no need for the compensating step \eqref{eq:scale-p_i+1}, since no subsequent polynomials will be affected.)

We can apply this scheme in a couple of alternative ways to accelerate the code:
\begin{itemize}
  \item Scale one of the $p_i$ by the reciprocal of its leading coefficient. This transforms it to a monic polynomial which, as observed in Section \ref{sec:monics}, eliminates one multiply instruction from the code (in the absence of FMA instructions). We can repeat this process for all but one of the polynomials, whose scaling factor is now fully constrained to be the reciprocal of the product of all the other scale factors.
  \item Scale each of the $p_i$ so that its leading coefficient has magnitude $u^{n_i}$, for a suitable constant $u$, with $n_i$ the degree of $p_i$. Now, where the unmodified algorithm would compute the value $x^a$ to be reused in each iteration, we can instead compute the product $u x^a$ at a cost of one extra multiply instruction and then reuse this value in every iteration to save one multiply per iteration.
\end{itemize}
In both cases, the saving in execution cost will be $m-1$ multiply instructions, but the computed result will be unaffected except by rounding differences.

%% file: implementation.tex
\section{Implementation details} \label{sec:implementation}

\subsection{Calculation of the ``magic constant"}

The FRSR and FRGR algorithms as presented contain a constant $c$, whose value is computed by Algorithm \ref{alg:FRGR-constants}. Recall that in section \ref{sec:standard_alg}, we alluded to a correspondence between $c$ and the so-called ``magic constant" of Listing \ref{lst:Quake} and the family of improved versions it has led to. Naming this latter constant $C$, we now develop a formula for computing $C$ from $c$. We shall assume a mantissa with $k_\text{mant}$ bits and an exponent bias of $k_\text{bias}$.

Corresponding to lines \ref{alg:FRGR:line:X} and \ref{alg:FRGR:line:y} of Algorithm \ref{alg:FRGR}, let $I(x)$ denote the value obtained when interpreting the bit pattern of the floating point number $x$ as an integer of the same word length, and let $F(X)$ denote the value obtained when the integer $X$ is interpreted as a floating point number.

We can use our pseudolog function $L(x)$ of equation \eqref{eq:L(x)} to write
\begin{align*}
  I(x) &= 2^{k_\text{mant}} (L(x) + k_\text{bias}) \,, \\
  F(X) &= L^{-1} (2^{-k_\text{mant}} X - k_\text{bias})\,.
\end{align*}
We can then transform lines \ref{alg:FRGR:line:X} - \ref{alg:FRGR:line:y} of  Algorithm \ref{alg:FRGR} into new versions which utilise the functions $I$ and $F$ and the constant $C$. Since the intermediate variables will now take on different values, we replace $X$ with $X'$, etc. The transformation can be written as
\begin{align*}
  X &= L(x)                                     & \xleftrightarrow{} && X' &= I(x) \,, \\
  Y &= \frac{c}{b} - \frac{a}{b} X & \xleftrightarrow{} && Y' &= C - \frac{a}{b} X' \,, \\
  y &= L^{-1} (Y)                           & \xleftrightarrow{} && y' &= F(Y') \,.
\end{align*}
Since we wish the value of the coarse approximation to remain unchanged by this transformation, we impose $y'=y$ and solve the resulting system for $C$ to obtain
\begin{equation} \label{eq:C}
  C = \frac{2^{k_\text{mant}}}{b} (c + k_\text{bias} (a+b))\,.
\end{equation}
If we now assume the IEEE 754 single-precision format, using the value $c=-\frac{1}{2}$ for FRSR obtained in Section \ref{sec:app-FRSR}, and applying equation \eqref{eq:C}, we find for $C$ a value $2^{21}\times761$, which can be written as the hexadecimal constant \texttt{0x5F200000}.

In practice, we cannot attain the theoretical accuracy quoted in that section, because evaluation error due to floating point rounding affects the result. However, there are some measures we can take to reduce the evaluation error, which are presented next.

\subsection{Implementation alternatives}

In this section we are interested in ways to reorder the operations without changing the execution cost. We will only consider the FRSR case, but the same technique can be applied to other cases.

Assuming that intermediate results are rounded to the nearest floating point value, there are often multiple ways to express a sequence of floating point operations all of which are algebraically equivalent but which may yield slightly different results depending on the order of operations. The IEEE 754 standard ensures that the result of the expression $f_0 * f_1$ will equal that of $f_1 * f_0$. However, the result of $(f_0 * f_1) * f_2$ may differ from that of $f_0 * (f_1 * f_2)$. In the FRSR case, this lack of associativity implies that we can change the order in which we compute the product $c_1 x y^2$ occurring in the evaluation of $p(z)$. We find that there are nine alternatives, which are easily enumerated:
\begin{align*}
  & c_1*x*y*y && c_1*y*x*y && c_1*y*y*x \\
  & x*y*c_1*y && x*y*y*c_1 && y*y*c_1*x \\ 
  & y*y*x*c_1 && (c_1*x)*(y*y) && (c_1*y)*(x*y)
\end{align*}
Another set of possible implementations can be derived by observing the way \citet{kadlec2010} factored the expression for the refinement step. By suitably transforming the coefficients $\{c_0,c_1\} \rightarrow \{c_0',c_1'\}$ and defining a temporary variable $w = c_1' * y$, we can use one of the following expressions for the refined approximation (Kadlec's version corresponding to the first expression):
\begin{align*}
  & w * (c_0' - x * y * y) && w * (c_0' - y * y * x) \\
  & w * (c_0' - x * y * w) && w * (c_0' - x * w * y) && w * (c_0' - y * w * x) \\
  & w * (c_0' - x * w * w) && w * (c_0' - w * w * x)
\end{align*}
We now have a total of 16 ways to implement the refinement computation, and some will yield better accuracy than others.\footnote{Note that some ways of ordering the products have implications for the domain on which the function retains its full accuracy. This will be addressed in the results section.}

Furthermore, in listing \ref{lst:Quake}, the right-shift operation is performed before the subtraction, but the order could be swapped to become
\[  Y = (C' - X) >> 1;  \]
where the modified magic constant $C'$ is whichever of $2 C$ or $2 C + 1$ yields the best result, effectively giving the constant one additional bit of precision.\footnote{The right-shift in this case should be performed as an unsigned operation since the topmost bit of $C'-X$ will often be $1$.}

\subsection{Alternative values for $c$}

We observed that Algorithm \ref{alg:FRGR-constants} gives us a free choice for the integer $s$ without affecting the value computed. In a practical application, though, this choice can affect the evaluation error.

Referring to Algorithm \ref{alg:FRGR}, if we increase $c$ by an amount $b$, $Y$ will increase by $1$, doubling $y$. This doubling is counteracted by scaling the coefficients of $p$ by appropriate powers of $2$. Since multiplication by a power of $2$ can be represented exactly, the final result is unaffected. However, adding any of $1, ..., b-1$ to $c$ modifies the mantissa of $y$ and there is in general no exact floating point representation for the correspondingly modified coefficients of $p$. Thus we have a set of $b$ alternative implementations (orthogonally with those of the previous section). In the FRSR case, we have two alternatives - $s$ even, and $s$ odd. As the results section will show, the best choice for the parity of $s$ in the classic FRSR is not the one which previous authors have traditionally chosen (although \citet{kadlec2010} alluded to the existence of such a choice).

\subsection{Tuning the values of $c, c_0, c_1$}

In a general setting, an analytic approach to minimising evaluation error by modifying the coefficients is a complex topic. The interested reader may refer to \citet{arzelier2019}, for example. One might consider a brute-force search over all possible combinations of constants, but even in the FRSR case the reader will quickly see that the search space is intractably large. Fortunately, it is possible to dramatically reduce the search space using considerations motivated by the methods of this paper (although we omit the details here). This approach can then be combined with the other alternatives explained in this section, and in doing so we find that the implementation yielding the best results is the function \texttt{FRSR\_Deg1()} shown in Section \ref{sec:results}, with a peak relative error of $6.501791 \times 10^{-4}$ (although see the note there regarding an alternative version).

%% file: results.tex
\section{Results} \label{sec:results}

Here we present several implementations for computing reciprocal square root with differing degrees for the refinement polynomial, together with a few examples for reciprocal and reciprocal cube root. All were found using the methods of this paper. The peak relative error is shown with each listing. Aside from the very first listing, for which the optimal solution had already been published, we believe all versions here yield greater accuracy than any others in their respective classes, demonstrating the usefulness of our techniques.

\subsection{Testing environment}

All tests were performed using code compiled with \texttt{gcc} for a 64-bit target, and run on an Intel Core i5 processor. This setup appears to be equivalent to that used by the majority of other authors, since the peak errors we find when testing their code agree precisely with the values published by those authors. Consistent with the work of other authors, the peak relative error of each function was measured over all positive normal floats (except where noted).

\subsection{Reciprocal square root}

A function returning the coarse approximation corresponds to a refinement step which uses the constant monic polynomial. In this case the optimal magic constant is the one found by \citet{lomont2003}:

\begin{minipage}{\linewidth}
\begin{lstlisting}[caption={$x^{-\frac{1}{2}}$ with degree-0 monic, $\epsilon$ =  $3.421284 \times 10^{-2}$}, label={lst:FRSR_Mon0}, language=C]
float FRSR_Mon0(float x)
{
  uint32_t X = *(uint32_t *)&x;
  uint32_t Y = 0x5F37642F - (X >> 1);
  return *(float *)&Y;
}
\end{lstlisting}
\end{minipage}

A general degree-0 polynomial incurs one extra multiply, but only yields a slight improvement in accuracy over the monic version:

\begin{minipage}{\linewidth}
\begin{lstlisting}[caption={$x^{-\frac{1}{2}}$ with degree-0 polynomial, $\epsilon$ = $2.943730 \times 10^{-2}$}, label={lst:FRSR_Deg0}, language=C]
float FRSR_Deg0(float x)
{
  uint32_t X = *(uint32_t *)&x;
  uint32_t Y = (0xBEBFFDAA - X) >> 1;
  float y = *(float *)&Y;
  return y * 0.79247999f;
}
\end{lstlisting}
\end{minipage}

With a degree-1 monic, the error reduces dramatically to $8.802292 \times 10^{-4}$, making it almost exactly twice as accurate as the original Quake code despite being faster than it by one multiply instruction:

\begin{minipage}{\linewidth}
\begin{lstlisting}[caption={$x^{-\frac{1}{2}}$ with degree-1 monic, $\epsilon$ = $8.802292 \times 10^{-4}$}, label={lst:FRSR_Mon1}, language=C]
float FRSR_Mon1(float x)
{
  uint32_t X = *(uint32_t *)&x;
  uint32_t Y = (0xBE167122 - X) >> 1;
  float y = *(float *)&Y;
  return y * (1.8909901f - x*y*y);
}
\end{lstlisting}
\end{minipage}

Using a general degree-1 polynomial corresponds to the classic FRSR algorithm. Our result improves on the \citet{kadlec2010} version:

\begin{minipage}{\linewidth}
\begin{lstlisting}[caption={$x^{-\frac{1}{2}}$ with degree-1 polynomial, $\epsilon$ = $6.501791 \times 10^{-4}$}, label={lst:FRSR_Deg1}, language=C]
float FRSR_Deg1(float x)
{
  uint32_t X = *(uint32_t *)&x;
  uint32_t Y = 0x5F5FFF00 - (X >> 1);
  float y = *(float *)&Y;
  return y * (1.1893165f - x*y*y*0.24889956f);
}
\end{lstlisting}
\end{minipage}

If we relax the requirement that same accuracy must hold all normal floats, there is a version even more accurate for $x < 1.8822997 \times 10^{38}$, beyond which it degrades slightly to a peak error of $ 6.502243 \times 10^{-4}$:

\begin{minipage}{\linewidth}
\begin{lstlisting}[caption={$x^{-\frac{1}{2}}$ with degree-1 polynomial, $\epsilon$ = $6.501686 \times 10^{-4}$ for $x < 1.8822997 \times 10^{38}$ }, label={lst:FRSR_Deg1Alt}, language=C]
float FRSR_Deg1Alt(float x)
{
  uint32_t X = *(uint32_t *)&x;
  uint32_t Y = 0x5F6004CC - (X >> 1);
  float y = *(float *)&Y;
  return y * (1.1891762f - y*y*x*0.24881148f);
}
\end{lstlisting}
\end{minipage}

A degree-$2$ monic uses one extra floating point add, but yields a $32$-fold improvement over our \texttt{FRSR\_Deg1()}, and an $86$-fold improvement over the original Quake code:

\begin{minipage}{\linewidth}
\begin{lstlisting}[caption={$x^{-\frac{1}{2}}$ with degree-2 monic, $\epsilon$ = $2.020644 \times 10^{-5}$}, label={lst:FRSR_Mon2}, language=C]
float FRSR_Mon2(float x)
{
  uint32_t X = *(uint32_t *)&x;
  uint32_t Y = 0x5F11107D - (X >> 1);
  float y = *(float *)&Y;
  float z = x*y*y;
  return y * (2.2825186f + z*(z-2.253305f));
}
\end{lstlisting}
\end{minipage}

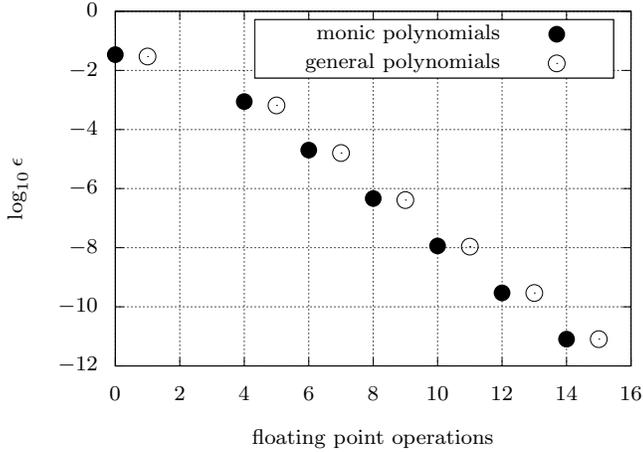
\begin{figure}
  \centering
  \input{monics}
  \caption{Comparison of peak relative error $\epsilon$ for monic vs. general polynomials used in the refinement step of the FRSR algorithm, assuming unlimited precision and no FMA. (The apparent gap in the data is the jump from a constant polynomial to a linear one, for which we must calculate $z$.)}
  \label{fig:monics}
\end{figure}

A general degree-2 polynomial (not shown here) only improves the peak error by roughly $25\%$ over a degree-2 monic. Motivated by this finding and the corresponding ones for degrees 0 and 1, we have computed the peak relative errors of FRSR versions using both monic and general polynomials up to degree 6, assuming unlimited precision. The comparison is shown in Figure \ref{fig:monics}, with the errors plotted against the number of floating point operations assuming no FMA is available. We find that the versions which employ monic polynomials are much better placed with respect to the cost/accuracy trade-off than those using general polynomials. Of particular note is degree 6, where the monic peak relative error $8.027828 \times 10^{-12}$ is highly comparable to that for a general polynomial, $8.027660 \times 10^{-12}$.

\subsection{Reciprocal}

Only the degree-1 version is shown, corresponding to the classic FRSR case:

\begin{minipage}{\linewidth}
\begin{lstlisting}[caption={$x^{-1}$ with degree-1 polynomial, $\epsilon$ = $1.116995 \times 10^{-4}$ for $x < 9.0209911 \times 10^{37}$}, label={lst:FRCP_Deg1}, language=C]
float FRCP_Deg1(float x)
{
  uint32_t X = *(uint32_t *)&x;
  uint32_t Y = 0x7FB504EC - X;
  float y = *(float *)&Y;
  return y * (0.6966215f - x*y*0.12130684f);
}
\end{lstlisting}
\end{minipage}

\subsection{Reciprocal cube root}

$x^{-\frac{1}{3}}$ is implemented here using a degree-1 general polynomial, corresponding to the classic FRSR case:

\begin{minipage}{\linewidth}
\begin{lstlisting}[caption={$x^{-\frac{1}{3}}$ with degree-1 polynomial, $\epsilon$ = $8.014543 \times 10^{-4}$}, label={lst:FRCR_Deg1}, language=C]
float FRCR_Deg1(float x)
{
  uint32_t X = *(uint32_t *)&x;
  uint32_t Y = 0x54638AFE - X/3;
  float y = *(float *)&Y;
  return y * (1.8696972f - (x*y)*(y*y)*1.2857759f);
}
\end{lstlisting}
\end{minipage}

Our degree-2 version improves on the one found by \citet{moroz2021} by using Algorithm \ref{alg:FRGR-constants} to derive an optimal value for the magic constant\footnote{It appears that \citet{moroz2021} chose a criterion which did not minimise $\rho$. Instead, by inverting equation \eqref{eq:C}, it is found that their magic constant corresponds to setting $t^*=t_1$; whereas, observing line \ref{alg:FRGR-constants:line:t_star} of Algorithm \ref{alg:FRGR-constants} we see that, for optimality, this value must be clamped to the interval [$\frac{1}{3}$,$\frac{2}{3}$]. Doing so yields $t^*=\frac{1}{3}$ and, together with $s=0$, the superior magic constant shown in listing \ref{lst:FRCR_Deg2}.}:

\begin{minipage}{\linewidth}
\begin{lstlisting}[caption={$x^{-\frac{1}{3}}$ with degree-2 polynomial, $\epsilon$ = $2.662789 \times 10^{-5}$}, label={lst:FRCR_Deg2}, language=C]
float FRCR_Deg2(float x)
{
  uint32_t X = *(uint32_t *)&x;
  uint32_t Y = 0x54B8E38E - X/3;
  float y = *(float *)&Y;
  float z = x*y*y*y;
  return y * (1.3739948f-z*(0.47285829f-z*0.092823250f));
}
\end{lstlisting}
\end{minipage}

The FRGR algorithm can also be used to implement $x^{-\frac{2}{3}}$. Only the version using a degree-1 general polynomial is presented:

\begin{minipage}{\linewidth}
\begin{lstlisting}[caption={$x^{-\frac{2}{3}}$ with degree-1 polynomial, $\epsilon$ = $1.190003 \times 10^{-3}$}, label={lst:FRCR2_Deg1}, language=C]
float FRCR2_Deg1(float x)
{
  uint32_t X = *(uint32_t *)&x;
  uint32_t Y = 0x69BC56FC - 2*X/3;
  float y = *(float *)&Y;
  float w = 0.8152238f * y;
  float v = x * w;
  return w * (1.7563311f - v*v*w);
}
\end{lstlisting}
\end{minipage}

A more accurate way to approximate $x^{-\frac{2}{3}}$ (though on a narrower domain) is to pass $x*x$ to the function \texttt{FRCR\_Deg1()}, at the cost of an additional multiply. (Squaring the \textit{result} of \texttt{FRCR\_Deg1()} is \textit{less} accurate.) Either method can then be used to calculate a cube root using $x * x^{-\frac{2}{3}}$. If using \texttt{FRCR2\_Deg1()} for this purpose, the final multiplication by $x$ can be eliminated by reworking the code.

\subsection{Multiple iterations for reciprocal square root}

Our 2-iteration version of FRSR is established by tuning the coefficients and instruction sequence for the first iteration, then finding the actual range $[z_{min}, z_{max}]$ of $z$ values produced by the resulting code, and then using this range to tune the second iteration in a separate pass. It would almost certainly be possible to obtain a slightly superior result by tuning all 5 numbers and both instruction sequences in concert.

\begin{minipage}{\linewidth}
\begin{lstlisting}[caption={$x^{-\frac{1}{2}}$ with 2nd iteration, $\epsilon$ =  $4.612440 \times 10^{-7}$}, label={lst:FRSR_Iter}, language=C]
float FRSR_Iter(float x)
{
  uint32_t X = *(uint32_t *)&x;
  uint32_t Y = 0x5F5FFF00 - (X >> 1);
  float y = *(float *)&Y;
  y *= 1.1893165f - x*y*y*0.24889956f;
  y *= 1.4999996f - (0.49999934f*y)*(x*y);
  return y;
}
\end{lstlisting}
\end{minipage}

We now apply the transformation of coefficients described in Section \ref{sec:acc-multi-iter} to make the second polynomial monic, eliminating one multiply instruction. The resulting code has execution cost equal to that of the original 2-iteration version of the Quake code, but more than 10 times better accuracy:

\begin{minipage}{\linewidth}
\begin{lstlisting}[caption={$x^{-\frac{1}{2}}$ with monic 2nd iteration, $\epsilon$ =  $4.639856 \times 10^{-7}$}, label={lst:FRSR_IterFast}, language=C]
float FRSR_IterFast(float x)
{
  uint32_t X = *(uint32_t *)&x;
  uint32_t Y = 0x5F5FFF00 - (X >> 1);
  float y = *(float *)&Y;
  y *= 0.9439607f - x*y*y*0.19755164f;
  y *= 1.8898820f - x*y*y;
  return y;
}
\end{lstlisting}
\end{minipage}

%% file: monics.tex
\begingroup
  \inputencoding{cp1252}%
  \fontfamily{Times-New-Roman}%
  \selectfont
  \makeatletter
  \providecommand\color[2][]{%
    \GenericError{(gnuplot) \space\space\space\@spaces}{%
      Package color not loaded in conjunction with
      terminal option `colourtext'%
    }{See the gnuplot documentation for explanation.%
    }{Either use 'blacktext' in gnuplot or load the package
      color.sty in LaTeX.}%
    \renewcommand\color[2][]{}%
  }%
  \providecommand\includegraphics[2][]{%
    \GenericError{(gnuplot) \space\space\space\@spaces}{%
      Package graphicx or graphics not loaded%
    }{See the gnuplot documentation for explanation.%
    }{The gnuplot epslatex terminal needs graphicx.sty or graphics.sty.}%
    \renewcommand\includegraphics[2][]{}%
  }%
  \providecommand\rotatebox[2]{#2}%
  \@ifundefined{ifGPcolor}{%
    \newif\ifGPcolor
    \GPcolorfalse
  }{}%
  \@ifundefined{ifGPblacktext}{%
    \newif\ifGPblacktext
    \GPblacktexttrue
  }{}%
  \let\gplgaddtomacro\g@addto@macro
  \gdef\gplbacktext{}%
  \gdef\gplfronttext{}%
  \makeatother
  \ifGPblacktext
    \def\colorrgb#1{}%
    \def\colorgray#1{}%
  \else
    \ifGPcolor
      \def\colorrgb#1{\color[rgb]{#1}}%
      \def\colorgray#1{\color[gray]{#1}}%
      \expandafter\def\csname LTw\endcsname{\color{white}}%
      \expandafter\def\csname LTb\endcsname{\color{black}}%
      \expandafter\def\csname LTa\endcsname{\color{black}}%
      \expandafter\def\csname LT0\endcsname{\color[rgb]{1,0,0}}%
      \expandafter\def\csname LT1\endcsname{\color[rgb]{0,1,0}}%
      \expandafter\def\csname LT2\endcsname{\color[rgb]{0,0,1}}%
      \expandafter\def\csname LT3\endcsname{\color[rgb]{1,0,1}}%
      \expandafter\def\csname LT4\endcsname{\color[rgb]{0,1,1}}%
      \expandafter\def\csname LT5\endcsname{\color[rgb]{1,1,0}}%
      \expandafter\def\csname LT6\endcsname{\color[rgb]{0,0,0}}%
      \expandafter\def\csname LT7\endcsname{\color[rgb]{1,0.3,0}}%
      \expandafter\def\csname LT8\endcsname{\color[rgb]{0.5,0.5,0.5}}%
    \else
      \def\colorrgb#1{\color{black}}%
      \def\colorgray#1{\color[gray]{#1}}%
      \expandafter\def\csname LTw\endcsname{\color{white}}%
      \expandafter\def\csname LTb\endcsname{\color{black}}%
      \expandafter\def\csname LTa\endcsname{\color{black}}%
      \expandafter\def\csname LT0\endcsname{\color{black}}%
      \expandafter\def\csname LT1\endcsname{\color{black}}%
      \expandafter\def\csname LT2\endcsname{\color{black}}%
      \expandafter\def\csname LT3\endcsname{\color{black}}%
      \expandafter\def\csname LT4\endcsname{\color{black}}%
      \expandafter\def\csname LT5\endcsname{\color{black}}%
      \expandafter\def\csname LT6\endcsname{\color{black}}%
      \expandafter\def\csname LT7\endcsname{\color{black}}%
      \expandafter\def\csname LT8\endcsname{\color{black}}%
    \fi
  \fi
    \setlength{\unitlength}{0.0500bp}%
    \ifx\gptboxheight\undefined%
      \newlength{\gptboxheight}%
      \newlength{\gptboxwidth}%
      \newsavebox{\gptboxtext}%
    \fi%
    \setlength{\fboxrule}{0.5pt}%
    \setlength{\fboxsep}{1pt}%
    \definecolor{tbcol}{rgb}{1,1,1}%
\begin{picture}(4320.00,3600.00)%
    \gplgaddtomacro\gplbacktext{%
    }%
    \gplgaddtomacro\gplfronttext{%
      \csname LTb\endcsname
      \put(-305,2041){\rotatebox{-270}{\makebox(0,0){\strut{}$\log_{10} \epsilon$}}}%
      \put(2375,154){\makebox(0,0){\strut{}$\textrm{floating point operations}$}}%
      \csname LTb\endcsname
      \put(3332,3206){\makebox(0,0)[r]{\strut{}$\textrm{monic polynomials}$}}%
      \csname LTb\endcsname
      \put(3332,2986){\makebox(0,0)[r]{\strut{}$\textrm{general polynomials}$}}%
      \csname LTb\endcsname
      \put(300,704){\makebox(0,0)[r]{\strut{}$-12$}}%
      \csname LTb\endcsname
      \put(300,1150){\makebox(0,0)[r]{\strut{}$-10$}}%
      \csname LTb\endcsname
      \put(300,1596){\makebox(0,0)[r]{\strut{}$-8$}}%
      \csname LTb\endcsname
      \put(300,2042){\makebox(0,0)[r]{\strut{}$-6$}}%
      \csname LTb\endcsname
      \put(300,2487){\makebox(0,0)[r]{\strut{}$-4$}}%
      \csname LTb\endcsname
      \put(300,2933){\makebox(0,0)[r]{\strut{}$-2$}}%
      \csname LTb\endcsname
      \put(300,3379){\makebox(0,0)[r]{\strut{}$0$}}%
      \csname LTb\endcsname
      \put(432,484){\makebox(0,0){\strut{}$0$}}%
      \csname LTb\endcsname
      \put(918,484){\makebox(0,0){\strut{}$2$}}%
      \csname LTb\endcsname
      \put(1404,484){\makebox(0,0){\strut{}$4$}}%
      \csname LTb\endcsname
      \put(1890,484){\makebox(0,0){\strut{}$6$}}%
      \csname LTb\endcsname
      \put(2376,484){\makebox(0,0){\strut{}$8$}}%
      \csname LTb\endcsname
      \put(2861,484){\makebox(0,0){\strut{}$10$}}%
      \csname LTb\endcsname
      \put(3347,484){\makebox(0,0){\strut{}$12$}}%
      \csname LTb\endcsname
      \put(3833,484){\makebox(0,0){\strut{}$14$}}%
      \csname LTb\endcsname
      \put(4319,484){\makebox(0,0){\strut{}$16$}}%
    }%
    \gplbacktext
    \put(0,0){\includegraphics[width={216.00bp},height={180.00bp}]{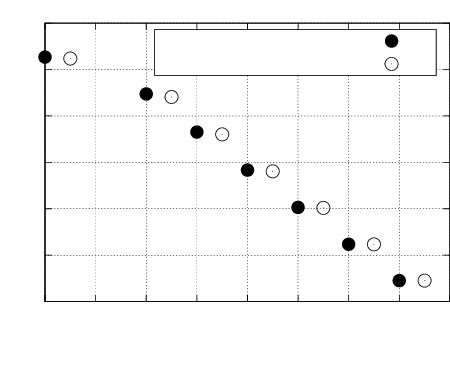}}%
    \gplfronttext
  \end{picture}%
\endgroup

%% file: acknowledgements.tex
\section*{Acknowledgements}

The author wishes to extend his sincere gratitude to Andreas Fredriksson, Mike Acton, and Jasmine Banks, without whose encouragement this project would certainly not have been completed.

%% file: appendices.tex
\begin{appendices}

\section{Proof that $\zeta$ is increasing with respect to $k$} \label{appendix:zeta}

Let $k \in \mathbb{Z^+}$ and $h>0$.
Using the binomial expansion
\[  \left( 1 + \frac{h}{k} \right)^k = \sum_{i=0}^{k} \binom{k}{i} \left( \frac{h}{k} \right) ^ i\,,  \]
we can re-express the term with index $i$ as
\[  \binom{k}{i} \left( \frac{h}{k} \right) ^ i = \frac{h^i}{i!}\frac{k(k-1)...(k-i+1)}{k^i} = \frac{h^i}{i!} \prod_{j=1}^{i-1} \left( 1-\frac{j}{k} \right)\,.   \]
For $i \in \{0,1\}$, the product on the right is empty and so this term is independent of $k$. But, for $i \geqslant 2$, it is strictly increasing in $k$. Furthermore, the number of such terms increases with $k$, and all terms are positive. Hence, $\left( 1 + \frac{h}{k} \right)^k$ is strictly increasing in $k$. Observing also that the expression always equals $1$ when $h=0$, we can say that $\left( 1 + \frac{h}{k} \right)^k$ is weakly increasing in $k$ when $k \in \mathbb{Z^+}$ and $h \geqslant 0$. Now replacing $h$ with $r+t$ and multiplying by $2^{s-r}$, we see from \eqref{eq:zeta} that $\zeta_{r,k}(c)$ is weakly increasing with respect to $k$, and that when at least one of $r$ and $t$ is non-zero, it is strictly increasing.

\section{Bounds on $\phi(k)$} \label{appendix:phi}

Suppose $x>0$. Expanding $e^x$ as a power series, we have
\begin{align} \label{eq:phi-bounds-aux}
  e^x (x-2) + x+2 & = \sum_{i=0}^{\infty} \frac{x^i}{i!} (x - 2) + x + 2 \nonumber \\
                              & = \sum_{i=3}^{\infty} \frac{(i-2)x^i}{i!} \nonumber \\
                              & > 0\,,
\end{align}
with the last line following since all terms in the sum are positive. Dividing by the positive number $2x(e^x-1)$ and rearranging, this becomes
\[  \frac{1}{e^x-1} - \frac{1}{x} + \frac{1}{2} > 0\,.  \]
On substituting $x = \frac{\ln{2}}{k}$, we find
\[  \frac{1}{2^\frac{1}{k}-1} - \frac{k}{\ln{2}} + \frac{1}{2} > 0 \,,  \]
which, recalling the definition of $\phi(k)$ (equation \eqref{eq:phi}), we can express as
\begin{equation} \label{eq:phi-lower-bound}
  \phi(k) > \left( \frac{1}{\ln{2}} - 1 \right) k + \frac{1}{2}\,.
\end{equation}
Noting that $k>0 \implies x>0$, the inequality \eqref{eq:phi-lower-bound} provides a lower bound for $\phi(k)$ valid for all positive $k$.

To find a suitable upper bound, we begin by defining a linear function
\[  \phi_0(k) = \left( \frac{1}{\ln{2}} - 1 \right) k + 2 - \frac{1}{\ln{2}} \,,  \]
and proceed to take the derivative of $\phi(k) -  \phi_0(k)$. We find it to be
\begin{equation} \label{eq:phi-phi1-deriv}
  \phi'(k) -  \phi_0'(k) = \frac{ 2^{\frac{1}{k}} \ln{2} } { k^2 \left( 2^{\frac{1}{k}}-1 \right) ^2 } - \frac{1}{\ln{2}}\,.
\end{equation}
We also have for $x>0$ that
\begin{align*}
  x^2 + 2 - 2 \cosh{x} & = x^2 + 2 - 2 \left( 1 + \frac{x^2}{2!} + \sum_{i=2}^{\infty} \frac{x^{2i}}{(2i)!} \right) \\
                                     & = -2 \sum_{i=2}^{\infty} \frac{x^{2i}}{(2i)!} \\
                                     & < 0 \,.
\end{align*}
Multiplication by $e^x$ yields
\[  x^2 e^x - (e^x - 1)^2 < 0\,,  \]
and on dividing by $(e^x-1)^2 \ln{2}$ we obtain
\[  \frac{x^2 e^x}{(e^x-1)^2 \ln{2}} - \frac{1}{\ln{2}} < 0\,.  \]
We now perform the same substitution as before, $x = \frac{\ln{2}}{k}$, and use \eqref{eq:phi-phi1-deriv} to find that for $k>0$ we have
\[  \phi'(k) - \phi_0'(k) < 0\,.  \]
Furthermore, $\phi(k) - \phi_0(k)$ is clearly continuous for $k>0$ and takes the value $0$ for $k=1$.

We conclude that for $k \geqslant 1$ we have $\phi(k)-\phi_0(k) \leqslant 0$, that is,
\begin{equation} \label{eq:phi-upper-bound}
  \phi(k) \leqslant \left( \frac{1}{\ln{2}} - 1 \right) k + 2 - \frac{1}{\ln{2}}\,,
\end{equation}
providing the desired upper bound. We summarise the results by combining \eqref{eq:phi-lower-bound} and \eqref{eq:phi-upper-bound}, giving
\begin{equation} \label{eq:phi-bounds}
  k \geqslant 1 \implies \left( \frac{1}{\ln{2}} - 1 \right) k + \frac{1}{2} < \phi(k) \leqslant \left( \frac{1}{\ln{2}} - 1 \right) k + 2 - \frac{1}{\ln{2}}\,.
\end{equation}

\section{Bounds on $t_0(k)$} \label{appendix:t0}

To find bounds on $t_0(k)$, we will make use of the following easily verified relationship between the functions $\phi$ and $t_0$ (see equations \eqref{eq:t0} and \eqref{eq:phi}):
\begin{equation} \label{eq:phi-t0-relationship}
  k>1 \implies \phi \left( \frac{k}{k-1} \right) = \frac{t_0(k)}{k-1} + 1\,.
\end{equation}
We begin by showing that the second derivative of $\phi(k)$ is positive for $k>0$. The second derivative is found to be
\begin{equation} \label{eq:phi-second-deriv}
  \phi''(k) = \frac { 2^\frac{1}{k} \ln{2} \left( 2k(1-2^\frac{1}{k}) + (2^\frac{1}{k}+1 )\ln{2} \right) } { k^4 (2^\frac{1}{k}-1)^3 }\,.
\end{equation}
Now given $x>0$, and making the substitution $x = \frac{\ln{2}}{k}$ in inequality \eqref{eq:phi-bounds-aux}, we obtain
\begin{equation} \label{eq:phi-2nd-deriv-aux}
  2 (1-2^\frac{1}{k}) + \frac{\ln{2}}{k} (2^\frac{1}{k}+1) > 0\,.
\end{equation}
Since $k>0$ precisely when $x>0$, the quantity $\frac{ 2^\frac{1}{k} \ln{2}  } { k^3 (2^\frac{1}{k}-1)^3 }$ will be positive, and if we multiply inequality \eqref{eq:phi-2nd-deriv-aux} by it, on comparing with \eqref{eq:phi-second-deriv} we see that we have now proven
\[  k>0 \implies \phi''(k)>0\,.  \]
For positive $k$, therefore, the curve of $\phi(k)$ is concave upward, and so any line which is tangent to it also bounds it from below. Consider the tangent at $k=1$, which we will denote as $\phi_1(k)$ and which can be written as
\begin{equation} \label{eq:phi1}
  \phi_1(k) = (2\ln{2}-1)k + 2(1-\ln{2})\,.
\end{equation}
We then have
\[  k>0 \implies \phi(k) - \phi_1(k) > 0\,.  \]
But now if $k>1$, then $\frac{k}{k-1} > 1$ also, hence we have
\[  \phi \left( \frac{k}{k-1} \right) - \phi_1 \left( \frac{k}{k-1} \right) > 0\,.  \]
Using \eqref{eq:phi-t0-relationship} and \eqref{eq:phi1}, this becomes
\[  \frac{t_0(k)}{k-1} + 1 - (2\ln{2}-1) \left( \frac{k}{k-1} \right) - 2(1-\ln{2}) > 0 \,.  \]
Multiplying through by $k-1$ and simplifying yields
\begin{equation} \label{eq:t0-lower-bound}
  t_0(k) > 2\ln{2}-1 \,,
\end{equation}
providing a lower bound for $t_0(k)$ valid for $k>1$.

By a similar token, since the curve of $\phi(k)$ is concave upward, any line which cuts the curve in exactly two places must lie above the curve on the interval between the two intersection points.

Using $\phi_2(k)$ to denote the line having intersection points at $k=1$ and $k=2$, we can express this line as
\begin{equation} \label{eq:phi2}
  \phi_2(k) = (\sqrt{2}-1) k + 2 - \sqrt{2}\,,
\end{equation}
allowing us to write
\[  1 \leqslant k \leqslant 2 \implies \phi(k) - \phi_2(k) \leqslant 0\,.  \]
Given $k \geqslant 2$, then we have $1 < \frac{k}{k-1} \leqslant 2$, and hence
\[  \phi \left( \frac{k}{k-1} \right) - \phi_2 \left( \frac{k}{k-1} \right) \leqslant 0\,.  \]
Using \eqref{eq:phi-t0-relationship} and \eqref{eq:phi2}, this becomes
\[  \frac{t_0(k)}{k-1} + 1 - (\sqrt{2}-1) \left( \frac{k}{k-1} \right) - 2 + \sqrt{2} \leqslant 0 \,.  \]
Now multiplying by $k-1$ and simplifying, we obtain the upper bound
\begin{equation} \label{eq:t0-upper-bound}
  t_0(k) \leqslant \sqrt{2} - 1 \,,
\end{equation}
valid for $k \geqslant 2$.

Finally, we combine \eqref{eq:t0-lower-bound} with \eqref{eq:t0-upper-bound} to obtain the desired result:
\begin{equation} \label{eq:t0-bounds}
  k \geqslant 2 \implies 2\ln{2}-1 < t_0(k) \leqslant \sqrt{2} - 1\,.
\end{equation}

\section{Upper bound on minimax error} \label{appendix:eps-bound}

Suppose a function $f$ is both continuous and positive on a given closed interval, and that $p$ is a minimax polynomial approximation (with respect to relative error) to $f$ on that interval. We will show that the resulting minimax error is less than $1$.

To see this, note that, whatever the degree of $p$, it is as least as accurate an approximation as any constant polynomial. Hence the result is proven if it can be shown to hold for the minimax polynomial of degree $0$.

Let $f_0$ and $f_1$ be the minimum and maximum values of $f$ on the given interval, and suppose we approximate $f$ by a constant $c_0$. Let $e_0$ and $e_1$ be the relative errors corresponding to $f_0$ and $f_1$, respectively. Then we have
\[  e_0 = 1-\frac{c_0}{f_0}\,, \quad e_1 = 1-\frac{c_0}{f_1}\,.  \]
If $c_0$ is the degree-$0$ minimax polynomial, equioscillation implies that $e_0 = -e_1$, whence we obtain
\[  c_0 = \frac{2 f_0 f_1}{f_0 + f_1} \,,  \]
and the minimax error is therefore
\[  |e_0| = |e_1| = \frac {f_1-f_0}{f_0+f_1} \,.  \]
Since $-f_0 < f_0$, then $f_1-f_0 < f_0+f_1$ and hence the minimax error must be strictly less than $1$.

\end{appendices}